\documentclass{conm-p-l}

\usepackage{amssymb,amsmath,amsthm,amscd,float,color,mathtools, pdflscape}
\allowdisplaybreaks

\usepackage{tabularx, pdflscape}

\usepackage{hyperref}
\usepackage{multicol}
\usepackage{algorithm}
\usepackage{algorithmic}
\usepackage{amsmath,amsbsy,amsfonts,amsopn,amstext} %, cite}
\usepackage{euscript,amssymb,amsbsy,amsfonts,amsthm,latexsym,amsopn,amstext,amsxtra,euscript,amscd}
\usepackage{graphicx, verbatim}

\usepackage{longtable}

\usepackage{graphicx}

\newtheorem{thm}{Theorem} 

\newtheorem{prop}{Proposition} 
\newtheorem{lem}{Lemma}

\newtheorem{exa}{Example} 

\newtheorem{alg}{Algorithm} 
\newtheorem{rem}{Remark}
\newtheorem{cor}{Corollary}

\newcommand{\ch}[2]
{\begin{bmatrix}
 #1 \\
 #2\\
\end{bmatrix}}

\newcommand{\chr}[4]
{\begin{bmatrix}
 #1 & #2\\
 #3 & #4
\end{bmatrix}}

\DeclareMathOperator\Aut{Aut}
\DeclareMathOperator\bAut{\overline{Aut}}

\DeclareMathOperator\gal{Gal }

\DeclareMathOperator\h{h}
\DeclareMathOperator\mH{\mathfrak h}

\newcommand\Z{\mathbb Z}
\newcommand\Q{\mathbb Q}
\newcommand\R{\mathbb R}
\newcommand\C{\mathbb C}

\def\P{\mathbb P}

\newcommand\M{{\mathcal M}}
\newcommand\X{\mathcal X}            %\newcommand\X{\mathfrak X}
\newcommand\B{\mathcal B}
\def\O{\mathcal O}
\def\L{\mathcal L}

\newcommand\m{\mathfrak m}

\def\p{\mathfrak p}
\newcommand\Jac{\mbox{Jac }}

\newcommand\iso{\cong}

\newcommand\emb{\hookrightarrow }

\newcommand\embd{\hookrightarrow}

\def\iso{\cong}

\def\>{\rangle}
\def\<{\langle}

\newcommand\D{\Delta}                % Discriminant
\newcommand\om{\omega}

\def\d{{\delta }}
\def\a{\alpha}

\def\t{\tau}

\def\k{\bar k}

\def\l{\lambda}

\usepackage{amssymb,amsmath,amsthm,amscd,float,color,mathtools, pdflscape}
\allowdisplaybreaks

\hypersetup{
  colorlinks   = true, %Colours links instead of ugly boxes
  urlcolor     = blue, %Colour for external hyperlinks
  linkcolor    = blue, %Colour of internal links
  citecolor   = red %Colour of citations
}

\usepackage[msc-links]{amsrefs}

\def\v{\mathfrak v}
\def\dis{\mathfrak D}

\def\L{\mathcal L}

\title{Rational points in the moduli space of genus two}

\author{L. Beshaj}
\address{Department of Mathematics \\
The University of Texas at Austin,
Austin, TX 78712
}
\email{beshaj@math.utexas.edu}

\author{R. Hidalgo}
\address{Departamento de Matem\'atica y Estad\'{\i}stica, Universidad de La Frontera, Temuco, Chile.}
\email{ruben.hidalgo@ufrontera.cl}

\author{S. Kruk}
\address{Department of Mathematics and Statistics \\
Oakland University \\
Mathematics and Science Center,  \\
Rochester, MI 48309}
\email{kruk@oakland.edu}

\author{A. Malmendier}
\address{
Department of Mathematics\\
Utah State University,
Logan, UT 84322
}
\email{andreas.malmendier@usu.edu}

\author{S. Quispe}
\address{Departamento de Matem\'atica y Estad\'{\i}stica, \\  Universidad de La Frontera, Temuco, Chile.}
\email{saul.quispe@ufrontera.cl}

\author{T. Shaska}
\address{Department of Mathematics and Statistics \\
Oakland University\\
Mathematics and Science Center,   \\
Rochester, MI 48309}
\email{shaska@oakland.edu}

\date{\today} % delete this line to display the current date

\def\p{\mathfrak p}
\def\dd{J_{16}}
\def\B{J_{30}}

\DeclareMathOperator\Gal{Gal }
\DeclareMathOperator\twist{Twist}
\DeclareMathOperator\homo{Hom}

\def\T{\th}
\def\th{{\theta}}

\def\bP{\mathbb P}

\def\mh{\mathfrak h}

\def\IA{\mathfrak A}
\def\IB{\mathfrak  B}
\def\IC{\mathfrak C}
\def\ID{\mathfrak  D}

\begin{document}

\begin{abstract}
We build a database of  genus 2  curves defined over $\Q$ which contains all curves with minimal absolute height $\h \leq 5$, all curves with moduli height $\mH \leq 20$, and all curves with extra automorphisms in standard form $y^2=f(x^2)$ defined over $\Q$ with height $\h \leq 101$. 
For each isomorphism class in the database, an equation over its minimal field of definition is provided, the automorphism group of the curve, Clebsch and Igusa invariants.  
The distribution of rational points in the moduli space $\M_2$ for which the field of moduli is a field of definition is discussed and some open problems are presented. 

\end{abstract}

\maketitle

\setcounter{tocdepth}{1}

%\tableofcontents

\def\x{\textbf{x}}

%**************************
\section{Introduction}

In \cite{MR3525576}  were  introduced concepts of \textit{minimum absolute height} of a binary form, the \textit{moduli height}, and discussed relations between the two.  Moreover, some computations were performed for binary sextics of minimum absolute height one.  A natural problem in that paper was to check whether a large database of binary forms (equivalently genus 2 curves) of relatively small minimum absolute height could be constructed. Comparing the minimum absolute height with the moduli height would be the main point of this database. Moreover, such a database would shed some light to other problems related to the rational points of the moduli space $\M_2$. For example, not every rational point $\p \in \M_2$ has a representative genus 2 curve $\X$ defined over $\Q$.  What is the percentage of points of bounded height for which the field of moduli is not a field of definition?  How does this ratio is affected when the height increases? 

From another point of view, we can list genus two curves based on their moduli height. For example, traditionally  there have been plenty of effort to count the curves with bounded discriminant. That is not an easy problem, but would it make more sense to have an estimate on the number of curves with bounded moduli height?  After all, the moduli height is the most natural way of sorting points in $\M_2$. 

And then, there are also the curves with automorphisms.  In \cite{Beshaj:2016ac} it was shown that such curves can always be written in an equation such  that the corresponding binary sextic is reduced; see \cite{MR3525574} for details. Reduced usually means minimal absolute height for the binary form. Is this really supported computationally? 

Our goal was to construct a database of genus 2 curves which addresses some of these questions.

\noindent In this paper we construct three main databases:  
\begin{itemize}
\item[i)] integral binary sextics of minimum absolute height $\leq 4$, 

\item[ii)] integral binary sextics with moduli height $\mh \leq 20$, 

\item[iii)] integral binary sextics  $f(x^2, z^2)$ of height $\h \leq 101$.  
\end{itemize}
All combined we have over 1 million isomorphism classes of genus two curves (equivalently rational points in $\M_2$) without counting twists. We compute the minimal absolute height, the moduli height, the discriminant, automorphism group, field of definition, and an equation over the field of definition for all such points. 
We discuss all the technical details of organizing the data and some open questions on the distribution of the rational points with non-trivial obstruction in the moduli space $\M_2$.

Let $k$ be an algebraically closed  field of characteristic zero and $\M_g$ the moduli space of smooth, projective, genus $g\geq 2$ algebraic curves defined over $k$.  
The moduli space $\M_2$ of genus 2 curves is the most understood moduli space among all moduli spaces.  This is mostly due to two main facts; first all genus two curves are hyperelliptic and therefore studying them it is easier than general curves, secondly even among hyperelliptic curves the curves of genus two have a special place since they correspond to  binary sextics which, from the computational point of view, are relatively well understood compared to higher degree binary forms.

Some of the main questions related to $\M_2$ have been to recover a nice equation for any  point $\p \in \M_2$.  Since $\M_2$ is a coarse moduli space, such equation is not always defined over the field of moduli of $\p$.   Can we find a universal equation for genus two curves over their minimal field of definition?  Can such equation provide a minimal model for the curve?  Does the height of this minimal model has any relation to the projective height of the corresponding moduli point $\p \in \M_2$?  What is the distribution in $\M_2$ of points $\p$ for which the field of moduli is not a field of definition?  The answers to these questions are still unknown.  

In \cite{data} we provide a computational package for computing with genus 2 curves and a  database of genus 2 curves which contains all curves  with height $\h \leq 5$, curves with moduli height $\mh \leq 20$, and curves with automorphism and height $\leq 101$. They are organized in three Python directories $\L_i$, $i=1, 2, 3$ as explained in Section~\ref{sec-11}.

The database is build with the idea of better understanding $\M_2$, the distribution of points in $\M_2$ based on the moduli height, the distribution of points for which the field of moduli is not a field of definition.  

The goal of this paper is twofold: to provide the mathematical background for most of the algorithms in \cite{data} and to discuss some of the open questions and problems raised there. 
  Most of the material of the first part has already appeared in the vast literature on genus two curves some of which is previous work of these authors.  For sake of completeness and straightening out some notational confusion we define all the basic invariants of genus two curves in this paper.

The current database and all the functions are implemented in Sage.  It improves and expands a previous Maple genus 2 computational algebra package as in \cite{arith-gen-2}.  There is another database of genus 2 curves in \cite{Booker:2016aa} which collects all genus 2 curves with discriminants $\leq 1000$.  Some remarks on how the two databases overlap can be found in the last section. 

There is a lot of confusion in the literature about the invariants of genus two curves.  We go to great lengths to make sure that all the invariants are defined explicitly and there is no room for misunderstanding.  We like to warn the reader that our invariants are different from the ones used by Magma and all the papers which use Magma in their computations. 

%*****************
%%\newpage

%Notice that throughout the paper a curve is identified by its moduli point $\p=(j_1, j_2, j_3)$.  This is somewhat different than previous works of these authors or the vast majority of literature on genus 2 curves, where other invariants are used.  The benefit of Igusa functions $j_1, j_2, j_3$ is that they are defined everywhere.  Their disadvantage is that, given as rational functions in terms of the coefficients of the curve,  they have larger degrees than other invariants (i.e. $i_1, i_2, i_3$ used in \cite{deg2} and other papers). 

%\bigskip

%\noindent \textbf{Notation:}   By a "curve" we always mean the isomorphism class of smooth, irreducible genus two curve $\X$  over $\C$.  This class is identified with the corresponding point $\left[  J_2 : J_4 : J_6 : J_{10} \right] $ in the weighted projective space $\mathbb{WP}^3_{(2,4,6,10)}$ or equivalently the affine point $\p=(i_1, i_2, i_3)$ (cf. Eq.~\eqref{mod-point}. All invariants of the curve can be computed by such moduli point.  Moreover, all invariants and functions mentioned in this paper are implemented in \cite{data}. 

%*************************************

%**********************
%\newpage
\section{A database of integral binary sextics}\label{sec-11}

Our main goal is to create an extensive database of integral binary sextics with  minimal absolute height and all the twists with minimal height. We will use the definitions of minimal absolute height, moduli height, and the basic properties of heights of polynomials from \cite{MR3525576} and will provide details in the next coming sections.   The database will be organized in  a Sage/Python dictionary, where the key will be the moduli point
\[
\p = 
\left\{
\aligned
& ( -1, i_1, i_2, i_3 ), \; \quad if \; J_2 \neq 0\\
& ( 0, t_1, t_2, t_3 ),  \; \qquad if \; J_2=0,  \\
%& \frac{J_{6}^5}{J_{10}^3} \;\; if \; J_2=0, J_4=0, J_6\neq 0  \\
%& \frac{J_{4}^5}{J_{10}^2} \;\; if \; J_2=0, J_6=0, J_4\neq 0\\
\endaligned
\right.
%\end{split}
\]
%
%First, let us see how we can  recover all the other invariants from this key. 
see Eq.~\eqref{mod-point}  for definitions of such moduli point.   The data is organized in   three   main dictionaries: 
\begin{itemize}
\item[i)]  integral binary sextics with minimum absolute height  $\h \leq 10$, 

\item[ii)] decomposable integral binary sextics  $f(x^2, z^2)$ with minimum absolute  $\h \leq 101$ and,

\item[iii)]  integral binary sextics with moduli height $\mh \leq 20$ \\
\end{itemize}
Each point in the database has the following invariants  \\
\[
\begin{aligned}
\p= \left(r , i_1, i_2, i_3 \right): \qquad        &  \h      & =  &\textit{ minimal absolute height }  \\
                              &  \mh         & =  & \textit{  moduli height }\\
                              &  \D           & =  & \textit{  minimal discriminant }\\
                              & \Aut (\p)    & = & \textit{ automorphism group }\\
                              &  C           & =  &  \textit{  Conductor }\\
                              &  M_{\p}           & = & \textit{  field of definition of the universal curve }\\
                              & \twist       & =   &          \textit{ List of  twists }  \\
\end{aligned}
\]
An entry in each dictionary  looks as the following:
\[ \left(r, i_1, i_2, i_3 \right): \left[ \h, \,  \mh,\,  \D, \, \Aut (\p) , \,  C, \, M_{\p}, \, \left[ [a_0, \dots , a_6], \dots , [b_0, \dots , b_6] \frac{} {} \right]   \frac {} {} \right]\]
We illustrate with an example.
\begin{exa}
Let $\X$ be  the curve with equation
\begin{equation}\label{example}
 y^2= x^6 -14x^4-82x^2+1
\end{equation} 
If we load the database in Sage and let 
\[ f= t^6 -14 t^4-82 t^2+1\]
then the command $\p$=\textbf{ModPoint (f)} displays
\[ \left( -1,   - \frac{49281147}{5410276}, \frac {706232480445}{12584301976}, \frac {3071021069999403}{17429644021121376256 }  \right) \]
If we ask if $\p \in \L$, where $\L$ is the second dictionary from above, the answer will be \textbf{yes} and $\L (\p)$ will display
\begin{small}
\[ 
\begin{split}
%& \left( -1,   - \frac{49281147}{5410276}, \frac {706232480445}{12584301976}, \frac {3071021069999403}{17429644021121376256 }  \right):\\
&  [82, 2^{14} \cdot  1163^5 , 17^2 \cdot 12301^2, [4, 2], C, \Q, [[1, 0, -14, 0, -82, 0, 1]]]
\end{split}
\]
\end{small}
which means that the minimal absolute height is $\h=82$, automorphism group with GapId $[4, 2]$ which is the Klein group $V_4$, minimal field of definition $\Q$, and minimal discriminant $\D = 17^2 \cdot 12301^2$. The moduli height is 
$\mh = 2^{14} \cdot  1163^5$. We chose not to display the conductor and all the twists. 
\end{exa}

%Notice that in this case the equation of the curve does have minimal absolute height as expected by results in \cite{beshaj-16}.

%For every curve given by an equation $y^2=f(x)$   we can display all its information by using the function \textit{Info(f)}.  In \ref{app-C} such information is displayed for the above curve. 

In the Appendix~\ref{app-B} is given a list of all the functions used for the genus 2 curves package in Sage. 
%The database contains all the genus two curves with height $\h < 10$, all the curves with moduli height $\mh < 21$, and all the genus 2 curves with extra automorphisms and height $\h \leq 101$.  
In the next few sections we will go over the necessary definitions and procedures to construct such databases.  The details will be explained in Section~\ref{data}.

%&&&&&&&&&&&&&&&&&&&&&&&&&&&&&&&&&&&&&&&&&&&&&&&&&&&&&&&&&&&&&&&&&&

%&&&&&&&&&&&&&&&
\section{Heights of genus two curves}

In this section we define heights on algebraic curves when such curves are given by some affine equation.   Throughout this paper $K$ denotes an algebraic number field and $\O_K$ its ring of integers.

Let $\X_g$ be an irreducible algebraic curve with affine equation $F(x, y)=0$ for $F(x, y) \in K [x, y]$.  We define the \textbf{height of the  curve over $K$} to be
\[H_K(\X_g):= \min \left \{    H_K(G) \, : H_K(G) \leq H_K(F) \right \}. \]
where the curve $G(x, y) =0$ is isomorphic to $\X_g$ over $K$.

If we consider the equivalence over $\bar K$ then we get another height which we denote it as $\overline H_K (\X_g)$ and call it \textbf{the height over the algebraic closure}. Namely,
\[ \overline H_K(\X_g)= \min \{H_K(G): H_K(G) \leq H_K(F)\},\]
 where the curve $G(x, y)=0$ is isomorphic to $\X_g$ over $\overline K$.

In the case that $K=\Q$ we do not write the subscript $K$ and use $H(\X_g)$ or  $\overline H(\X_g)$.  Obviously, for any algebraic curve $\X_g$ we have $\overline H_K(\X_g) \leq H_K(\X_g)$. In \cite{MR3525576}  is proved that  given $K$ a number field such that $[K:\Q] = d$, the height $H_K(\X_g)$  and   $\overline H_K(\X_g)$    are  well defined.

\begin{thm}[\cite{MR3525576}]\label{Thm1}
Let $K$ be a number field such that $[K:\Q] \leq d$. Given a constant $\h_0 \geq 1$  there are only finitely many curves 
%(up to isomorphism over $\bar K$)  
such that $H_K(\X_g) \leq \h_0$.
\end{thm}

As an immediate corollary we have the following

\begin{cor}
Let $\h_0 \geq 1$ be a fixed integer,  $K$ a number field,  and  $\O_K$ its ring of integers.  For any  genus $g \geq 2$ curve $\X_g$ defined over $\O_K$ with height $\h (\X_g) = \h_0$ there are only finitely many twists of $\X_g$ with height $\h_0$. 
\end{cor}

%We can say a lot more about superelliptic curves and in particular about genus 2 curves in Weierstrass form. 

%\begin{thm}
%Let $\X$ be a genus 2 curve given by an equation in Weierstrass form. Then,

%i)  The number of curves of height $\h_0$ is ...

%ii)  The number of curves of height $\h_0$ defined over $\Q$ is ....

%iii) The number of twists of height $\h_0$ is ...

%\end{thm}

%************************************
%%\newpage
%\subsection{Computing  the height $H(\X_g)$ of a genus $g\geq 2$ curve $\X_g$.}

Given a genus two curve $\X_g$ the following algorithm computes  a curve isomorphic over $K$ to $\X_g$ of minimum height

\begin{alg}

\textbf{Input:} an algebraic curve $\X_g : F(x, y)=0$, where 
 $F$ has degree $d$ and is  defined over $K$

\textbf{Output:} an algebraic curve $\X_g^\prime : G(x,y)=0$ such that
$\X_g^\prime \iso_K \X_g$ and $\X_g^\prime$ has minimum height.\\

\textbf{Step 1:} Compute $c_0 =H_K(F)$

\textbf{Step 2:} List all points $P \in \P^{s}(K)$ such that $H_K(P) \leq c_0$.

\textbf{Note:}  $s$ is the number of terms of $F$ which is the number of monomials
of degree $d$ in $n$ variables, and this is equal to $\binom{d+n-1}{d}$. From  Thm.~\ref{Thm1}  there are only finitely many such points assume $P_1, \dots, P_r$.

\textbf{Step 3:} for $i=1$ to $r$ do

\hspace{15mm}  Let $G_i(x, y) = p_i$;

\hspace{20mm} if $g(  G_i(x, y) ) = g(\X_g)$ then

\hspace{24mm}  if $G_i(x, y) =0\iso_K F(x,y)=0$

\hspace{28mm}   then add $G_i$ to the list $L$

\hspace{24mm} end if;

\hspace{20mm}end if;

\textbf{Step 4: } Return all entries of $L$ of minimum height , $L$ has curves isomorphic over $K$ to $\X_g$ of minimum height.

\end{alg}

\medskip

Note that this algorithm is not very efficient if we start with an algebraic curve of genus two and very big height.  Hence, the question that can be raised at this point is: how can we reduce the height of the curve? This is done using reduction theory, see \cites{MR3525574, SC} and others for more details,  and some elementary ways of reducing are given next.     The following  elementary lemma is  useful.

%Next we show some computational results for curves of height $ h < 10$.
%\subsubsection{Some elementary ways of reducing}

%There are some partial results when reducing the height of a curve.  The following  elementary lemma is  useful. 

\begin{lem}
Let $\X$ be a superelliptic curve  with Weierstrass equation $ y^m = \sum_{i=0}^d \, a_i x^i$, defined over $\Z$,  
and height $\h (\X)$. Let $p$ be a prime  such that $p \, | \, a_0$ and $\v_p (a_i ) = \alpha_i$,  such that $\alpha_0 \geq \alpha_i$, for $i=0, \dots , d$. Choose $m$ to be the largest nonnegative  integer which satisfies 
\[ m \leq \frac {\alpha_0 - \alpha_i} i, \quad i=1, \dots , d. \]
Then,  
%the image $\X^\sigma$ of the transformation  \; $ \sigma: (x, y)  \mapsto \left( p^m \cdot x, \,  y \right) $, 
%
there is a twist $\X^\prime$ of $\X$   such that  
\[ \v_p \left( \h(\X^\prime) \right) \leq  \v_p \left( \h(\X)   \right) - m.\]
\end{lem}

\proof  Let $p$ be a   prime such that
\[ a_i = p^{\a_i} \cdot b_i, \; \; \textit{such that} \; \; ( p, b_i) =1 \]
In other words, $\v_p (a_i ) = \alpha_i$, for $i=1, \dots , d$,  as in the assumptions of the theorem. 
%If $p^{\alpha_0} | a_i$, for $i=1, \dots , 6$, then we just divide by $p^{\alpha_0}$.  
Hence, the equation of the curve is 
\[ y^m=\sum_{i=0}^d  p^{\alpha_i} \, b_i \cdot  x^i \]
Choose $m$  as the largest  nonnegative integer such that 
\[ m  \leq  \frac {\alpha_0 - \alpha_i} i, \quad i=1, \dots , d. \]
If $m=0$, then the curve can no further reduced by this method at the prime $p$.   If  $m>0$, then we let
\[ (x, y)  \mapsto \left( p^m \cdot x, \, y \right) \]
which gives the curve  
\[ y^m= \sum_{i=0}^d  \;  a_i^\prime \, x^i       =  \sum_{i=0}^d  \;  p^{\alpha_i + im } \, b_i \,  x^i \]
Then,
\[ \v_p (a_i^\prime )= \alpha_i + i\cdot m. \]
Hence, to have reduction of $\v_p \left(\h (\X) \right)$ we must have $\alpha_i + i m \leq \alpha_0$ for $i=1, \dots , d$.  Thus, $m \leq   \frac {\alpha_0 - \alpha_i} i$, for $i=1, \dots d$.  Choosing the largest such $m$ will result to the biggest possible reduction on $p$. Dividing both sides by the content of the polynomial, which has the maximum power of $p$ as a factor,  gives a twist $\X^\prime$ of $\X$ with height which has valuation at $p$,  $\v_p \left(   \h(\X^\prime )   \right)    \leq \v_p \left(   \h (\X) \right) - m$. 
This completes the proof.
\qed

\begin{rem}
Notice that the curve $\X^\sigma$ could be in the same $\Gamma$-orbit of $\X$ or could be a twist of $\X$, depending on the values of $m$, where $\Gamma$ is the modular group.   
\end{rem}

\begin{cor}
Let $\X$ be a curve with $\Aut (\X)\iso D_6$ and equation 
\[ y^2 = x^6 + x^3 + s\]
where $s\in \Z$ such that it has a prime factorization 
\[ s= p^\alpha \cdot s^\prime, \quad \textit{where } \quad (s, s^\prime ) =1 \]
Then, we can reduce the  height by the transformation $x \mapsto p^m \cdot x$, where $m = \lfloor{ \frac \alpha 6   }  \rfloor$
\end{cor}

\begin{exa}
Let us consider the curve from Ex.~\ref{ex-d6}, namely
\[ y^2= x^6+ x^3 + 2^{33} \]
This curve has height $\h = 2^{33}$. Then $m$ has to be the largest nonnegative integer such that it is $\leq $ to 
\[ \frac {33- 0} 6, \frac {33-0} 5, \]
which makes $m=5$.

Consider the transformation $x\mapsto 2^5 \cdot x$.  Then the curve becomes
\[ y^2= 2^{15} \cdot x^6 + x^3 + 2^{18}\]
which is with height $\h = 2^{18}$.  
\end{exa}

Next we will define the moduli height of genus $g$ curves. 
%**************************************************
%%\newpage
\subsection{Moduli height of curves}
%
%In this section we define the height in the moduli space of curves and investigate how this height can be used to study the curves. Our main goal is to investigate if the height of the moduli point has any relation to the height of the curve.

Let $g$ be an integer $g \geq 2$ and $\M_g$ denote the coarse moduli space of smooth, irreducible algebraic curves of genus $g$. It is known that $\M_g$ is a quasi projective variety of dimension $3g-3$.  Hence, $\M_g$ is embedded in $\P^{3g-2}$. Let $\p \in \M_g$. We call the moduli height $\mH (\p)$ the usual height $H(P)$ in the projective space $\P^{3g-2}$.  Obviously, $\mH (\p)$ is an invariant of the curve.
In \cite{MR3525576} is proved the following result.

\begin{thm}For any constant $c\geq 1$, degree $d\geq 1$, and genus $g\geq 2$  there are finitely many superelliptic curves $\X_g$ defined over the ring of integers $\O_K$ of an algebraic number field $K$ such that   $[K:\Q] \leq d$ and  $\mH (\X_g) \leq c$.
\end{thm}

While the above theorem shows that the number of curves with bounded moduli height is finite, determine this number seems to be a very difficult problem. 

%\noindent \textbf{Open problem}  Determine the number of genus $g\geq 1$ curves with bounded moduli height $\lambda_0$.   

%**********************************************************
%%\newpage
\section{Genus 2 curves over \texorpdfstring{$\C$}{C}}

In this section we give a quick overview of the basic setup for genus two curves.  The material is part of the folklore on the literature of genus 2 curves and we don't mention all the possible references. While the main definitions and results on what follows are valid for any $g \geq 2$ we only state them for the case $g=2$. We mainly follow the approach of \cites{Cl, Bo, Ig, Ig-62, Ig-67, Me}.

\subsection{Periods and invariants}

Let $\X$ be a genus $g = 2$ algebraic curve. We choose a symplectic homology basis for $\X$, say $ \{ A_1,
 A_2, B_1,   B_2\},$ such that the intersection products $A_i \cdot A_j = B_i \cdot B_j =0$ and
$A_i \cdot B_j= \d_{i j},$ where $\d_{i j}$ is the Kronecker delta. We choose a basis $\{ w_i\}$ for the
space of holomorphic 1-forms such that $\int_{A_i} w_j = \d_{i j}$. The matrix 
\[   \Omega= \left[ \int_{B_i} w_j   \right] \]
is  the \emph{period matrix} of $\X$.  The columns of the matrix $\left[ I \ | \Omega \right]$ form a lattice $\Lambda$ in  $\C^g$ and the Jacobian  of $\X$ is $\Jac (\X) = \C^g/ \Lambda$.  Let $\mathbb H_g$ be the \emph{Siegel upper-half space}  and  $Sp_{4}(\Z)$ is the \emph{symplectic group}. Then $\Omega \in \mathbb H_g$.
\begin{prop}
Two period matrices $\Omega, \Omega^\prime$ define isomorphic principally polarized abelian varieties if and only if they are in the same orbit under the action of $Sp_{4}(\Z)$ on $\mathbb H_g$.
\end{prop}
Hence,  there is an injection
\[ \M_2 \embd \mathbb H_2/ Sp_{4}(\Z) =: \mathcal A_2 \] 
For any $z \in \C^2$ and $\t \in \mathbb H_2$ \emph{Riemann's theta function} is defined as
\[ \T (z , \t) = \sum_{u\in \Z^2} e^{\pi i ( u^t \t u + 2 u^t z )  } \]
where $u$ and $z$ are $2-$dimensional column vectors and the products involved in the formula are matrix products. The fact that the imaginary part of $\t$ is positive makes the series absolutely convergent over any compact sets. Therefore, the function is analytic.  The theta function is holomorphic on $\C^2\times {\mathbb H}_2$ and satisfies
\[\T(z+\tau u)=\T(z,\tau),\quad \T(z+\tau,\tau u)=e^{-\pi i( u^t \tau u+2z^t u )}\cdot   \T(z,\tau),\]
where $u\in \Z^2$.  Any point $e \in \Jac (\X)$ can be written uniquely as
$e = (b,a) \begin{pmatrix} 1_2 \\ \Omega \end{pmatrix}$, where $a, b \in \R^2$  are row vectors.  
We shall use the notation $[e] = \ch{a}{b}$ for the characteristic of $e$. For any $a, b \in \Q^2$, the theta function with rational characteristics is defined as
\[ \T  \ch{a}{b} (z , \t) = \sum_{u\in \Z^2} e^{\pi i ( (u+a)^t \t (u+a) + 2 (u+a)^t (z+b) )  }. \]
When the entries of column vectors $a^t$ and $b^t$ are from the set $\{ 0,\frac{1}{2}\}$, then the characteristics $ \ch {a}{b} $ are called the \emph{half-integer characteristics}. The corresponding theta functions with rational characteristics are called \emph{theta characteristics}. A scalar obtained by evaluating a theta characteristic at $z=0$ is called a \emph{theta constant}. 
%Points of order $n$ on $\Jac (\X)$ are called the $\frac 1 n$-\emph{periods}. 
%
Any half-integer characteristic is given by
\[
\m = \frac{1}{2}m = \frac{1}{2}
\begin{pmatrix} m_1 & m_2   \\ m_1^{\prime} & m_2^{\prime}    \end{pmatrix}
\]
where $m_i, m_i^{\prime} \in \Z.$  For $\gamma = \ch{\gamma ^\prime}{\gamma^{\prime \prime}} \in \frac{1}{2}\Z^{4}/\Z^{4}$ we define 
\[ e_*(\gamma) = (-1)^{4 (\gamma ^\prime)^t \gamma^{\prime \prime}}.\]
Then,
\[ \T [\gamma] (-z , \t) = e_* (\gamma) \T [\gamma] (z , \t).\]
We say that $\gamma$ is an  \textbf{even}  (resp.  \textbf{odd}) characteristic if $e_*(\gamma) = 1$ (resp. $e_*(\gamma) = -1$).
%;  see \cite{arith-gen-2} for details.

%&&&&&&&&&&&&&&&&&&&&&&&&&&&&&&&&&&&&&&&&&&&&&&&&&&&&&&&&&&&&&&&&&&&  theta functions and jacobians

For any genus 2 curve we have six odd theta characteristics and ten even theta characteristics. The following
are the sixteen theta characteristics, where the first ten are even and the last six are odd. For simplicity,
we denote them by $\theta_i = \ch{a} {b}$ instead of $\theta_i \ch{a} {b} (z , \t)$ where $i=1,\dots ,10$ for the
even theta functions.
\begin{small}
\[
\begin{split}
\theta_1 = \chr {0}{0}{0}{0} , \,  \theta_2 = \chr {0}{0}{\frac12} {\frac12} ,  \theta_3 =\chr
{0}{0}{\frac12}{0} , \, \, \theta_4 = \chr {0}{0}{0}{\frac12} , \, \,    \theta_5 = \chr{\frac12}{0}
{0}{0} ,\\
  \theta_6  = \chr {\frac12}{0}{0}{\frac12} , \, \,
  \theta_7 = \chr{0}{\frac12} {0}{0} , \, \,
  \theta_8 = \chr{\frac12}{\frac12} {0}{0} , \, \,
  \theta_9 = \chr{0}{\frac12} {\frac12}{0} , \, \,
  \theta_{10} = \chr{\frac12}{\frac12} {\frac12}{\frac12} ,\\
\end{split}
\]
\end{small}
and the odd theta functions  correspond to the following characteristics
\[   \chr{0}{\frac12} {0}{\frac12} , \,
   \chr{0}{\frac12} {\frac12}{\frac12} , \,
    \chr{\frac12}{0} {\frac12}{0} , \, \,
    \chr{\frac12}{\frac12} {\frac12}{0} , \,
    \chr{\frac12}{0} {\frac12}{\frac12} , \,
    \chr{\frac12}{\frac12} {0}{\frac12}  
\]
The complete set of thetanulls above are not independent, their relations are given via Frobenious relations. There are four theta constants which generate all the others, 
namely  \textbf{fundamental theta constants} $\theta_1, \, \theta_2, \, \theta_3, \, \theta_4$; see \cite{theta-1} for details. 
% They generate all the other theta constants. 
The following is Igusa's result, which is valid for any $g\geq 2$.  We only state it for $g=2$.
\begin{thm}
The complete set of theta constants uniquely determine the isomorphism class of a principally polarized abelian variety of dimension 2. 
\end{thm}
For curves of genus 2 this can be made more precise.   Let a genus 2 curve in Rosenheim form  be given by
\begin{equation} \label{Rosen2}
Y^2=X(X-1)(X-\lambda)(X-\mu)(X-\nu).
\end{equation} 
By the sa called Picard's lemma $\lambda, \mu, \nu$   can be written as follows:
\begin{equation}\label{Picard}
\l = \frac{\theta_1^2\theta_3^2}{\theta_2^2\theta_4^2}, \quad \mu = \frac{\theta_3^2\theta_8^2}{\theta_4^2\theta_{10}^2}, \quad \nu =
\frac{\theta_1^2\theta_8^2}{\theta_2^2\theta_{10}^2}.
\end{equation}
We can determine an equation of the curve in terms of the fundamental thetas as follows:
\begin{prop}[\cite{theta-1}]
\label{possibleCurve}
Every genus two curve can be written in the form:
\begin{equation}
y^2 = x \, (x-1) \, \left(x - \frac {\theta_1^2 \theta_3^2} {\theta_2^2  \theta_4^2}\right)\, \left(x^2 \, -   \frac{\theta_2^2 \, \theta_3^2 +
\theta_1^2 \, \theta_4^2} { \theta_2^2 \, \theta_4^2} \cdot    \a  \, x + \frac {\theta_1^2 \theta_3^2} {\theta_2^2 \theta_4^2} \, \a^2 \right),
\end{equation}
where $\a = \frac {\theta_8^2} {\theta_{10}^2}$ and in terms of $\, \, \theta_1, \dots , \theta_4$ is given by
\[
  \a^2 + \frac {\theta_1^4 + \theta_2^4 - \theta_3^4 - \theta_4^4}{\theta_1^2 \theta_2^2 - \theta_3^2 \theta_4^2 } \, \a + 1 =0
\]
Furthermore,  if $\alpha = {\pm} 1$ then $\X$ has an extra involution. 
\end{prop}
From the above we have that $\theta_8^4=\theta_{10}^4$ implies that $\X$ has an extra involution.  Hence, the Klein viergrouppe  $V_4 \embd Aut(\X)$.   The last part of the lemma above shows that if $\theta_8^4=\theta_{10}^4$ then all coefficients of the genus 2  curve are given as rational functions of the 4 fundamental theta functions. Such fundamental theta functions determine the field of moduli of the given curve. Hence, the curve is defined over its field of moduli.

\begin{cor}
Let $\X$ be a genus 2 curve which has an extra involution. Then $\X$ is defined over its field of moduli.
\end{cor}
We will revisit the curves defined over their field of moduli again in the coming sections. 

%*******************************************
\subsection{Siegel modular forms}

Here we define Siegel modular forms 
$ \psi_4$, $\psi_6$, $\chi_{10}$, $\chi_{12} $, of degree  4, 6, 10, and 12, 
as in \cite{Ig-67}*{pg. 848}. 
\begin{equation}\label{Siegel}
\begin{split}
2^2\cdot \psi_4 & = \sum (\theta_m)^8 \\
2^2 \cdot \psi_6 & = \sum_{syzygous} \pm \left( \theta_{m_1}  \theta_{m_2}   \theta_{m_3} \right)^4 \\
-2^{14} \cdot \chi_{10} & = \prod (\theta_m)^2  \\
2^{17} \cdot 3 \cdot \chi_{12} & = \sum  \left( \theta_{m_1}  \theta_{m_2}  \cdots  \theta_{m_6} \right)^4 \\
2^{39}\cdot 5^3 \cdot \chi_{35} & = \left(  \prod \theta_m\right )  \left(\sum_{azygous} \pm \left( \theta_{m_1}  \theta_{m_2}   \theta_{m_3} \right)^{20}\right).  \\
\end{split}
\end{equation}
In definition of $\chi_{12}$ the summation is taken over all G\"opel systems as explained in \cite{theta-1}, where all the G\"opel systems are displayed.

%The following functions 
%
%\[
%\begin{split}
% j_1 (\tau)  & = 2 \cdot 3^5 \,  \frac {\chi_{12}^5} {\chi_{10}^5}\\
% j_2 (\tau) &   = 2^{-3} \cdot 3^3 \, \frac {\psi_4 \chi_{12}^3} {\chi_{10}^4}    \\
% j_3 (\tau) & = 2^{-5} \cdot 3 \, \left( \frac {\psi_6 \chi{12}^2} {\chi_{10}^3}   + 2^2 \cdot 3 \frac {\psi \chi_{12}^3 } {\chi_{10}^4  } \right) \\
% \end{split}
% \]
%
% are defined in  \cite[pg. 4]{Lauter-Yang}. 

Theta constants provide a complete system of invariants for isomorphism classes of principally polarized varieties of dimension $g=2$.  But there are two main issues with this approach:
First, these invariants are not independent.  This can be fixed via the fundamental theta constants as in Prop.~\ref{possibleCurve}, however computationally things get difficult when we try to express all the results in terms of $\theta_1, \theta_2, \theta_3, \theta_4$. 
%Relations how fundamental theta constants can be expressed in terms of other functions can be found in \cite{theta-1, MR3525572}.
%
Secondly, and more importantly, they are defined analytically. Naturally, we would like to have algebraically defined invariants.  

%***************************
%%\newpage
\section{Algebraic invariants}

% While this is considered basic folklore by the experts, there is plenty of confusion in the literature about the set of invariants used to identify the isomorphism classes of genus 2 curves. Several sets of invariants are used starting with Bolza, Clebcsh, Igusa and continuing in the last two decades with Wamelen, Shaska, Gaudry,   and many many papers that have followed.
%For details on genus 2 curves we refer the reader to \cite{sh-13}.

Let $f(x, z)$ be a binary sextic defined over a field $k$,  $\mbox{char } k =0$,  given by
\begin{equation}\label{eq_1}
\begin{split}
f(x,z)       &=      a_0x^6+a_1x^5z+ \dots +a_6z^6       \\
&      =  (z_1x-x_1z)(z_2x-x_2z) \dots (z_6x-x_6z)
\end{split}
\end{equation}
A \textbf{covariant} $I$ of $f(x,z)$ is a homogenous polynomial in $x,z$
with coefficients in $k[a_0, \dots , a_{2g+2}]$. The \textbf{order } of $I$ is the degree of $I$ as a polynomial in
$x,z$ and the  \textbf{degree  }of $I$ is the degree of $I$ as a polynomial in $k[a_0, \dots , a_{2g+2}]$. An \textbf{invariant} is a covariant of  order zero. The binary form $f(x,z)$ is a covariant of order $2g+2$ and degree 1.
%; see \cite{vishi} for details.  
Throughout this paper we will use as basic references \cite{Cl}, \cite{Bo}, and  \cites{Ig, Ig-62, Ig-67}.

%**************
\subsection{Invariants and covariants via transvections}
%******************************************************

For any two binary forms $f$ and $g$ the symbol $(f, g)_r$ denotes the $r$-transvection.
%; see \cite{vishi} for details. 
Notice that the transvections are conveniently computed in terms of the coefficients of the binary forms.  

Let $f(x, z)$ be a binary sextic as in Eq.~\eqref{eq_1} and consider the  following covariants
\begin{equation}\label{Y-transv}
\begin{aligned}
& \Delta  = \left( (f, f)_4, (f, f)_4  \right)_2,     &  Y_1  = \left(f,  (f, f)_4  \right)_4 \\
& Y_2  = \left( (f, f)_4, Y_1  \right)_2,    &  Y_3  = \left( (f, f)_4, Y_2  \right)_2  \\
\end{aligned}
\end{equation}
The \textbf{Clebsch invariants} $A, B, C, D$  are defined as follows
\begin{equation}
\begin{aligned}
 & A  = (f, f)_6, \;                                & B  = \left( (f, f)_4,  (f, f)_4  \right)_4, \\
 & C  = \left( (f, f)_4, \Delta \right)_4, \;    & D  = \left(  Y_3, Y_1 \right)_2   \\
\end{aligned}
\end{equation}
see Clebsch \cite{Cl} or Bolza \cite{Bo}*{Eq.~(7), (8), pg. 51}  for details. 
%
%We also have the following covariants
%
%\begin{equation}
%x  = (Y_2 Y_3)_1, \quad y  =  (Y_3Y_1)_1, \quad z  = (Y_1Y_2)_1 
%\end{equation}
%
%where $Y_1, Y_2$ and $Y_3$ are as in Eq.~\eqref{Y-transv}.    

We display the invariants $A, B, C, D$ in terms of the coefficients in the Appendix. The following result is elementary but very important in our computations.

%******************************
%%\newpage
\subsection{Root differences}
%*********

Let $f(x, z)$ be a binary sextic as above and set
$D_{ij}:=   \begin{pmatrix} x_i & x_j \\  z_i & z_j \end{pmatrix}$.
For $\tau  \in SL_2(k)$,  we have
\[ \tau (f) =  (z_1^{'} x  - x_1^{'} z)   \dots  (z_6^{'}  x - x_6^{'}  z), \quad  \textit{ with } \quad
 \begin{pmatrix}   x_i^{'}  \\ z_i^{'}  \end{pmatrix}  = \tau^{-1} \, \begin{pmatrix} x_i\\ z_i \end{pmatrix}.
\]
Clearly $D_{ij}$ is  invariant under this action of $SL_2(k)$ on $\mathbb  P^1$. Let $\{i, j,  k, l, m, n  \} =\{ 1, 2, 3,  4, 5, 6 \}$. Treating $a_i$ as  variables, we construct the following elements in the ring of invariants $\mathcal R_6$
%; see \cite{Ig} .  
%In   \cite{vishi} a much simpler proof is provided   using geometric invariant theory for all $\mbox{char} k \neq 2, 3$. 
%Define the following invariants 
%
\begin{small}
\begin{equation}\label{j-invariants}
\begin{split}
\IA  & =  a_0^2 \,  \prod_{fifteen}   (12)^2 (34)^2 (56)^2  = \displaystyle \sum_{\substack {i<j,k<l,m<n}}  D_{ij}^2D_{kl}^2D_{mn}^2 \\
& \\
\IB  & =   a_0^4 \, \prod_{ten}  (12)^2 (23)^2 (31)^2 (45)^2 (56)^2 (64)^2 = \sum_{\substack {i<j,j<k, \\ l<m,m<n}}    D_{ij}^2D_{jk}^2D_{ki}^2D_{lm}^2D_{mn}^2D_{nl}^2 \\
 & \\
\IC  & =  a_0^6 \,  \prod_{sixty}   (12)^2 (23)^2 (31)^2 (45)^2 (56)^2 (64)^2  (14)^2 (25)^2 (36)^2  \\
 & = \sum_{ \substack {i<j,j<k, l<m, m<n \\  i<l' , j < m', k<n' \\  l',m',n' \in \{ l, m, n\} } }
 D_{ij}^2 D_{jk}^2 D_{ki}^2 D_{lm}^2 D_{mn}^2 D_{nl}^2  D_{{il}^{'}}^2 D_{{jm}^{'}}^2 D_{{kn}^{'}}^2 \\
 & \\
\ID  & =  a_0^{10} \prod_{i<j} (i j)^2   \\
\end{split}
\end{equation}
\end{small}
These    invariants, sometimes called \textbf{integral invariants},   are defined in \cite{Ig}*{pg. 620} where they are denoted by $A, B, C, D$.  Incidentally even Clebsch invariants which are defined next are also denoted by $A, B, C, D$ by many authors.  

%\begin{lem} $\IA, \IB, \IC, \ID$ are $SL_2 (k)$ invariants and expressed in terms of the coefficients they are  given by:
%\begin{small}
%\begin{equation}
%\begin{split}
%\IA  = & - 2^3 \cdot 3 \cdot 5 \cdot  A \\ 
%\IB  = & 2 \cdot 3^2\cdot 5 \,    \left( 75  B  - 8  \cdot  A^2 \right)\\
%\IC  =  & 2^2 \cdot  3^3 \cdot 5 \,   ( 16 A^3 -  200 A B +   375  C )\\
%\ID =  & 2\cdot 3^4   (6000   A^3B  - 374A^5 + 1060 A^2   C - 18750 AB^2   - 38500 BC - 28125 D )\\
%\end{split}
%\end{equation}
%\end{small}
%\end{lem}

%\proof  Is this correct?   Do we have a reference for this?
%\qed

%The following invariants 
%
%\[
%I_2   =  \IA, \quad 
%I_4 = (4I_{2}^2-\IB), \quad 
%I_6 = (8I_{2}^3-160I_2I_4-\IC), \quad
%I_{10}  =   \D \\
%
%\]

To quote Igusa \textit{"if we restrict to integral invariants, the discussion will break down in characteristic 2 simply because Weierstrass points behave badly under reduction modulo 2"}; see \cite{Ig}*{pg. 621}.  Next we define invarints which will work in every characteristic. 

%************************

\subsection{Igusa invariants} % $J_2, J_4, J_6, J_{10}$}

In \cite{Ig}*{pg. 622} Igusa defined what he called \textbf{basic arithmetic invariants}, which are now commonly known as  \textbf{Igusa invariants}
\[ 
\begin{aligned}
J_2 &= \frac 1 {2^3} \IA, \;                                      & J_4=  \frac 1 {2^5 \cdot 3} (4J_2^2-\IB) , \\
J_6 &= \frac 1 {2^6 \cdot 3^2} (8J_2^3-160J_2 J_4 -\IC), \;       &   J_{10}= \frac 1 {2^{12}} \ID  \\
\end{aligned}
\]
While most of the current literature on genus 2 curves uses invariants $\IA, \IB, \IC, \ID$, which are now most commonly labeled as $I_2, I_4, I_6, I_{10}$, Igusa went to great lengths in \cite{Ig} to define $J_2, J_4, J_6, J_{10}$ and to show that they also work in characteristic 2. 
\begin{lem}\label{lem3}
$J_{2i}$ are homogeneous elements in $\mathcal R_6$ of degree $2i$,  for $i$ = 1,2,3,5.
\end{lem}
A degree $d \geq 2$ binary form $f(x, z)$  is called \textbf{semistable} if it has no root of multiplicity $> \frac d 2$.

%\begin{lem}\label{lem4}
%A sextic has a root of multiplicity exactly three if and only if the basic invariants take the form
%\begin{equation}\label{eq6}
%I_2 = 3r^{2}, \quad I_4 =  81r^{4}, \quad I_6 = r^{6},  \quad I_{10} = 0.
%\end{equation}
%for some $ r \neq 0 $.
%\end{lem}

\begin{lem}\label{lem5}
A sextic has a root of multiplicity at least four if and only if the basic invariants vanish simultaneously.
\end{lem}

So a sextic for which all the basic invariants vanish simultaneously is not semistable.

Throughout this paper,  we will use these invariants  $J_2, J_4, J_6, J_{10}$. For  a binary form $f(x, z)=\sum_{i=0}^6 a_i \, x^i z^{6-i}$ as in Eq.~\eqref{eq_1}  we display   these invariants in terms of the coefficients $a_0, \dots , a_6$ in  \ref{app-A}.

Invariants $\, \, \{ J_{2i} \}$   are homogeneous polynomials of degree $2i$ in  $k[a_0, \dots , a_6]$, for  $i=1, 2, 3, 5$.   These $J_{2i}$    are invariant under the natural symmetries of the roots, namely the action of $S_6$ on the set of the roots.  They are also natural on the action of  $SL_2(k)$ on the sextic. 

\begin{lem}\cite{Ig}*{Cor. pg. 632}
Two genus 2 curves $C$ and $C^\prime$ are isomorphic if and only if there exists an $r\neq 0$ such that 
\[ J_{2i} (C) = r^{2i} \cdot J_{2i} ( C^\prime), \qquad \textit{ for } \quad i=1, 2, 3, 5.\]
\end{lem}

\begin{rem}
There is a lot of confusion in the literature about the definitions of the generators of $\mathcal R_6$.  Many authors have used different names and different notations for generators of $\mathcal R_6$. Moreover, many times the symbols $\IA, \IB, \IC, \ID$ are defined by scaling with different constants by many authors.  Hence,  sometimes equations in terms of these invariants might not be the same.  To avoid any confusion, we display the expressions of $j_2, J_4, J_6, j_{10}$ in terms of the coefficients of the binary sextic in Eq.~\eqref{eq_J}. 
\end{rem}

Igusa invariants are expressed in terms of the Clebsch invariants as follows:
\begin{equation}\label{J-Clebsch}
\begin{split}
 J_2  = &\; - 2^3 \cdot 3 \cdot 5  \, A\;, \\
 J_4  = & \;  2^3 \cdot 2^5 \, (75B-8A^2) \\
 J_6  = & \;  2^2 \cdot 3^3 \cdot 5 \; \left( 16 \, A^3 - 200 \, A \, B + 375 \, C \right)\\
 J_{10} = & \;   2^3 \cdot 4 \, \left(  - 384 \, A^5 + 6000 \, A^3 \, B + 10000 \, A^2 \, C  - 18750 \, A \, B^2 \right. \\
  & \; \; \; \qquad \left.  - 37500 \, B \, C - 28125 \, D  \right) \\  
\end{split}
\end{equation} 
Conversely, the invariants $(A,B,C,D)$ are polynomial expressions in the Igusa invariants $(J_2, J_4, J_6, J_{10})$ with rational coefficients are displayed in Eq.~\ref{iClebsch_invariants} in the Appendix.

%******* Ketu 

Thus, the invariants of a sextic define a point in a weighted projective space $[J_2 : J_4 : J_6 : J_{10}] \in \mathbb{WP}^3_{(2,4,6,10)}$. It was shown in \cite{vishi} that points in the projective  variety $\operatorname{Proj}\, \C [J_2, J_4, J_6, J_{10}]$ which are not on $J_2=0$ form the variety $\mathcal{U}_6 $ of moduli of sextics.
Equivalently, points in the weighted projective space $\{[J_2 : J_4 : J_6 : J_{10}] \in \mathbb{WP}^3_{(2,4,6,10)}: J_{10} \not = 0\}$ are in one-to-one correspondence with isomorphism classes of sextics.  

%***********

\subsection{Absolute invariants}

Dividing any $SL_2 (k)$  invariant by another one of the same degree gives an invariant under $GL_2(k)$ action.  
%The  term    \textbf{absolute invariants} is used first by Igusa \cite{Ig}   for $GL_2 (k)$ invariants.  

%Other sets of absolute invariants are used by different authors.  For example Gaudry/Schost \cite{gaudry} and Shaska et al.  \cites{deg2, deg3, ants} in their early papers have used 
%
We follow \cite{Ig-62}*{pg. 181} and define the following absolute invariants
\begin{equation}
i_1:=144 \frac {J_4} {J_2^2}, \quad i_2:=- 1728 \frac {J_2J_4-3J_6} {J_2^3}, \quad i_3 :=486 \frac {J_{10}} {J_2^5}
\end{equation}
for $J_2 \neq 0$. Notice that Igusa denotes them by $x_1, x_2, x_3$.

\begin{thm}\cite{Ig-62}*{Theorem 3}
The three absolute invariants can be expressed by the four modular forms in the form  
\begin{equation}\label{i-modular}
i_1 = \frac {\psi_4 \chi_{10}^2} {\chi_{12}^2}, \quad i_2 = \frac {\psi_6 \, \chi_{10}^3} {\chi_{12}^3}, \quad i_3= \frac {\chi_{10}^6 } {\chi_{12}^5} 
\end{equation}
\end{thm}
Hence, the theorem says that every modular form, in the degree 2 case, is expressed in the Eisenstein series of weight four, six, ten, and twelve. It is the main result of \cite{Ig-62}.

There is another set of absolute invariants 
\begin{equation} \label{abs_inv}
j_1: =  \frac{\IA^5}{\ID}, \quad j_2: = \frac{ \IB \IA^3 }{\ID} , \quad j_3 : =  \frac{\IC \IA^2}{\ID }  ,
\end{equation}
The popularity of invariants $j_1, j_2, j_3$ is due to P. Van Wamelen who used them in \cite{wamelen}  and later implemented them in Magma. They have no advantages over the invariants $i_1, i_2, i_3$ since both sets of invariants are not defined for $J_2 =0$ or equivalently $\IA=0$. Indeed, we could not find a direct proof that such invariants generate the field of invariants for $\mathcal R_6$, even though it is probably true. 
%
%Notice that there are some computational benefits in using $i_1, i_2, i_3$ when doing symbolic computations, since their degrees are smaller than those of $j_1, j_2, j_3$.   
The results of other authors who have worked with $j_1, j_2, j_3$ can be converted into our results and vice-versa using the formulas
\[ i_1 = 144 \, \frac {j_2} {j_1}, \quad i_2= - 1728 \, \frac {j_2 - 3 j_3} {j_1}, \quad i_3 = 486 \frac 1 {j_1} \]
and conversely
\[ j_1 = 486 \, \frac 1 {i_3}, \quad j_2 = \frac {27} {8} \, \frac {i_1} {i_3}, \quad j_3 = \frac 3 {32} \, \frac {i_2 + 12 i_1} {i_3} \]

\begin{rem}
 It is to be noted that for computational purposes  $i_1, i_2, i_3$ are much better than $j_1, j_2, j_3$ since they are of lower degrees.  Especially, for our purposes it is very important that the moduli point $\p =( r, i_1, i_2, i_3)$ (cf. next section) it is expressed in as small numbers as possible. 
\end{rem}

There is another set of invariants defined by Igusa in \cite{Ig},
\[ t_1= \frac {J_2^5} {J_{10}},\quad  t_2 = \frac {J_4^5} {J_{10}^2},  \quad t_3 = \frac {J_6^5} {J_{10}^3}\]
which are defined everywhere in the moduli space.  Due to their high degrees, they become difficult to use especially when someone  wants to do symbolic computations. 
%They were used in some computations in \cite{arith-gen-2}, but their high degrees make them unsuitable especially in symbolic computations. 

In the case $J_2=0$ we define
\begin{equation}\label{eq_a_1_a_2}
a_1:= \frac {J_4 \cdot J_6} {J_{10}}, \quad a_2:=\frac {J_6 \cdot J_{10}} {J_4^4}
\end{equation}
to determine genus two fields with $J_2=0$, $J_4\neq 0$, and $J_{6}\neq 0$ up to isomorphism.    They were used in \cites{Me, gaudry, deg2,   ants} and others. Moreover, when 
$J_2=0, J_4=0, J_6\neq 0 $ we have $\frac{J_{6}^5}{J_{10}^3}$ as an invariant and when $J_2=0, J_6=0, J_4\neq 0$ we have $\frac{J_{4}^5}{J_{10}^2}$ as an invariant.

%\begin{rem}
All our computations are made in terms of absolute invariants $i_1, i_2, i_3$   or the Igusa invariants $J_2, J_4, J_6, J_{10}$. In the Appendix, we provide formulas how to convert back and forth among all sets of invariants.
%\end{rem}
When  $J_2=0$ we will use invariants $t_1, t_2, t_3$ instead, in this case $t_1=0$, so this locus is determined by $t_2$ and $t_3$. 

%**********************************************************************

\subsection{Representing a point in the moduli space}

In creating a large database of genus two curves we need a way to identify uniquely a point $\p\in \M_2$.  Such point would ideally be defined everywhere, computationally feasable, and should not involve large integers.  
Of course, the representation $\p = (t_1, t_2, t_3)$ is unique, but these invariants are very large and extremely difficult to handle in a large database. The set of invariants $(i_1, i_2, i_3)$ is nicer, but they are not defined for $J_2=0$.

To find a compromise and simplify the implementation in \cite{data}, for  a given genus 2 curve $\X$ we define the corresponding \textbf{moduli point} $\p = [\X]$ to be
\begin{equation}\label{mod-point}
%\begin{split}
\p = 
\left\{
\aligned
& (i_1, i_2, i_3) \; \; if \; J_2 \neq 0\\
& (t_1, t_2, t_3) \; \; if \; J_2=0  \\
%& \frac{J_{6}^5}{J_{10}^3} \;\; if \; J_2=0, J_4=0, J_6\neq 0  \\
%& \frac{J_{4}^5}{J_{10}^2} \;\; if \; J_2=0, J_6=0, J_4\neq 0\\
\endaligned
\right.
%\end{split}
\end{equation}
This moduli point $\p$ together with $J_2$ will be the \textit{key} in our dictionary of genus 2 curves.  It will be stored as a 4-tuple $(r, \p )$ as explained before.  

%A projective moduli point will be the point  
%
%\[ \mathfrak P = \left[ J_2^5, 144 J_4 \, J_2^3, - 1728 (J_2\, J_4 - 3 J_6 ) \, J_2^2, 486 J_{10}\right] \]
%
%Notice that for a given  moduli point $\p = [ 1: i_1 : i_2 : i_3]$ we can always choose 
%
%\[ J_2 =1, \; J_4 = \frac {i_1} {144}, \; J_6= \frac 1 3 \left( \frac {i_2} {1728} + \frac {i_1} {144} \right), \; J_{10} = \frac {i_3} {486} \]

For $I_2 \not =0$ we use the variables $i_1, i_2, i_3$  to write
\begin{equation*}
\label{weighted_J}
\begin{split}
 \Big\lbrack J_2 : J_4 : J_6 : J_{10} \Big\rbrack =   \left\lbrack 1 : \frac{1}{2^4 \,3^2} \,  i_1 : \frac{1}{2^6 \, 3^4} \, i_2 + \frac{1}{2^4 \, 3^3} \, i_1
 :  \frac{1}{2 \cdot 3^5} \, i_3\right\rbrack \in \mathbb{WP}^3_{(2,4,6,10)} \;,
\end{split}
\end{equation*}
see \cite{MS-1} *{pg. 13} for details.  Notice that $i_1, i_2, i_3$ are denoted by $\mathbf{x}_1, \mathbf{x}_2, \mathbf{x}_3$ and $J_2, J_4, J_6, J_{10}$ respectively by $I_2, I_4, I_6, I_{10}$ in \cite{MS-1}.

%************************************************************

\subsection{Relations among different sets of invariants}
There are three sets of $SL_2(k)$ invariants most commonly used, namely $(A, B, C, D)$, $(\IA, \IB, \IC, \ID)$, and $(J_2, J_4, J_6, J_{10})$.  Many authors will use the notation $(I_2, I_4, I_6, I_{10})$ for $(\IA, \IB, \IC, \ID)$ due to the Magma notation (they are called  \textit{Igusa Clebsch invariants} in Magma) or some scaling of them.

Notice that  
\begin{equation} 
J_{2i}= \frac 1 {2^{4i}} \; I_{2i}, \qquad i=1, 2, 3, 5.
\end{equation}

For a given genus two curve $\X$ with equation $y^2=f(x)$, the relations between invariants $\IA, \IB, \IC, \ID$ and even Siegel modular forms $\psi_4, \psi_6, \chi_{10}, \chi_{12}, \chi_{35}$ of $\mathcal{A}_2$ are  given by 
Igusa in \cite{Ig-67}*{p.~\!848}:
\begin{equation}
\label{invariants}
\begin{aligned}
 \IA(f) & = -2^3 \cdot 3 \, \dfrac{\chi_{12}(\tau)}{\chi_{10}(\tau)} \;,                                            & \IB(f)  = \phantom{-} 2^2 \, \psi_4(\tau) \;,\\
 \IC(f) & = -\frac{2^3}3 \, \psi_6(\tau) - 2^5 \,  \dfrac{\psi_4(\tau) \, \chi_{12}(\tau)}{\chi_{10}(\tau)} \;,     & \ID(f)   = -2^{14} \, \chi_{10}(\tau) \;.
\end{aligned}
\end{equation}
In the Appendix are given the relations between $A, B, C, D$ and $J_2, J_4, J_6, J_{10}$. Using such relations we get the  absolute invariants  $i_1, i_2, i_3$ in terms of $A, B, C, D$ are as follows:
\begingroup\makeatletter\def\f@size{8}
\begin{align*}
i_1 & = \frac 9 {10} \cdot \frac {75B - 8A^2} {A^2} \\
i_2 & = \frac {27} {50} \cdot \frac { 112 A^3-900AB - 1125 C} {A^3} \\
i_3 & = \frac {81} {2^{13} \cdot 5^5} \cdot \frac { 384 A^5 - 6000 A^3B - 10 000 A^2C + 18 750 A B^2 + 37 500 B C + 28125 C D} {A^5} \\
\end{align*}
\endgroup
%%
%
%\[
%\begin{split}
%x & = \frac {2^3} {3^2 \cdot 5^2} \cdot {\frac {j_{{1}}+ 20\,j_{{2}}}{j_{{1}}}}\\
%y  & =  \frac {2^4} {3^3 \cdot 5^3} \cdot                    {\frac {j_{{1}}+2^4\cdot 5 \,j_{{2}}-2^3\cdot 3\cdot5^2 \,j_{{3}}}{j_{{1}}}}\\
%z & =  \frac {   2^6 } {3^4 \cdot 5^5} \cdot     
%\frac {9j_1^2 + 700j_1j_2 - 12400j_2^2 - 3600j_1j_3 + 48000j_2j_3 + 10800000j_1} {j_1^2}    
%\end{split}
%\]
%
%\begin{rem}
%While these formulas are quite useful, they are not defined everywhere. 
%Notice that it is exactly here that the implementation of the Mestre's algorithm fails for certain curves with $J_2=0$ (cf. Example~\ref{Me-fail}). 
%\end{rem}

Next we see how to construct the equation of a curve starting from a moduli point $\p \in \M_2$. 

%**********************************************
\subsection{Recovering an equation of the curve from the moduli point}
Some other invariants are as follows
\[ \begin{split}
&A_{ij} = (Y_iY_j)_2, \qquad \qquad \qquad  \quad (1 \leq i,j \leq 3)\\
&H_{ijk} =(fY_i)_2 (fY_j)_2 (fY_k)_2  \qquad (1 \leq i,j, k \leq 3)\\
&R = - (Y_1 Y_2)(Y_2Y_3)(Y_3Y_1)
\end{split}
\]
Note that $A_{ij}$ and  $H_{ijk}$  can be expressed in terms of the  Clebsch invariants, see \cite{Me}*{pg. 318}. 
Moreover, $R^2 = \frac 1 2 \det M$ where $M$ is the \textbf{Clebsch matrix} defined as follows. 
%************
%
\begin{equation}\label{M-matrix}
M= \begin{bmatrix}
A_{11}  &  A_{12}  &  A_{13}\\
A_{12} &  A_{22} & A_{23}\\
A_{13} & A_{23} & A_{33}
\end{bmatrix}\end{equation}
The following is an immediate consequence of classical invariant theory

Let $\B$ denote the degree 30 invariant 
\[ \B := \det M \]
An expression of $\B$ in terms of $J_2, J_4, J_6, J_{10}$ is given in \cite{deg2}*{Theorem 3}, where it is also expressed in terms of the roots of the sextic and its coefficients.

Let us now compute the determinant of the minor $ S= \begin{pmatrix} A_{1,1} & A_{1, 2} \\ A_{2, 1} & A_{2, 2} \end{pmatrix}$ from the matrix $M=\left[ A_{i, j} \right]$.
\[
\begin{split}    
\det S & = A_{1, 1} A_{2, 2} - A_{1, 2}^2 \\
& =  \frac {2^{16}}{ 3^6 \cdot 5^4 J_2^8}  \left(  15\, J_2^{3} J_4 \,  J_6  -4\, J_2^4 J_4^2   -175\, J_2^{2} J_4^{3} +2430\, J_{10}\, J_2^{3}-9\, J_2^{2} J_6^{2} \right.  \\
&  \left.+  1488\, J_2 \, J_4^{2} J_6 -64\, J_4^{4}+113400\, J_{10} \, J_2 \, J_4 -2880 \, J_4 \, J_6^{2} - 648000\, J_{10} \, J_6 \right)
\end{split}
\]
We define the new invariant 
\[ 
\begin{split}
\dd := & 15\, J_2^{3} J_4 \,  J_6  -4\, J_2^4 J_4^2   -175\, J_2^{2} J_4^{3} +2430\, J_{10}\, J_2^{3}-9\, J_2^{2} J_6^{2}  +  1488\, J_2 \, J_4^{2} J_6 \\
& -64\, J_4^{4}+113400\, J_{10} \, J_2 \, J_4 -2880 \, J_4 \, J_6^{2} - 648000\, J_{10} \, J_6\\
\end{split}
\]
which is a degree 16 invariant.

%**************************************************

\section{Automorphisms}
%*****************************************************
%
Let $\X$ be a genus 2 curve defined over an algebraically closed field $k$. Let $K$ be the function field of $\X$. Then $K$ has exactly one genus 0 subfield of degree 2, denote it $k(x)$. Since $k(x)$ is the only genus 0 subfield of degree 2 of $K$, then $G = \Aut(K/k)$ (or equivalently $\Aut (\X)$)  fixes $k(x)$. It is the fixed field of the hyperelliptic involution $\tau$ in $\Aut (K)$. Thus $\tau$ is in the center of $\Aut (K)$ and $\<\tau\> \lhd \Aut (K)$. 
The quotient group $\overline{\Aut} (K) =\Aut (K) / \<\tau \> $  is called the \textbf{reduced automorphism group} of $\X$ (or equivalently of $K$).

The reduced automorphism group embeds as a finite subgroup of $PGL(2, \C)$, therefore it is isomorphic to $C_n, D_n, A_4, S_4, A_5$.  Hence, the full automorphism group $\Aut (\X)$ is a degree 2 central extension of $\overline \Aut(\X)$.

The extension $K/k(x)$ is ramified at exactly six distinct points, namely the points $P = \{\om_1, \cdots, \om_6\}$ in $\bP^1$. The corresponding places of this points in $K$ are called Weierstrass points of $K$. The group $\Aut (K)$ permutes the 6 Weierstrass points and the group $\overline{\Aut} (K)$ permutes accordingly $\{\om_1, \cdots, \om_6\}$ in its action on $\bP^1$ as subgroup of $PLG(2, \C)$. This yields an embedding $\overline{\Aut} (K) \hookrightarrow S_6$, so all elements of $\overline{\Aut} (K)$ have order $\leq 6$.

In any characteristic different  from 2, the automorphism group $\Aut(\X)$ is isomorphic to one of the groups   given by the following theorem from \cite{deg2}.

\begin{thm}[\cite{deg2}]
\label{thm1}
The  automorphism group $G$ of a genus 2 curve $\X$ in characteristic $\ne2$ is isomorphic to \ $\Z_2$,
$\Z_{10}$, $V_4$, $D_8$, $D_{12}$, $SL_2 (3)$, $ GL_2(3)$, or $2^+S_5$. The case when $G \iso 2^+S_5$ occurs
only in characteristic 5. If $G \iso SL_2 (3)$ (resp., $ GL_2(3)$) then $\X$ has equation $Y^2=X^6-1$ (resp.,
$Y^2=X(X^4-1)$). If $G \iso \Z_{10}$ then $\X$ has equation $Y^2=X^6-X$.
\end{thm}

An \textbf{elliptic involution} of $K$ is  an  involution in $G$ which is different from $z_0$ (the hyperelliptic involution). Thus the elliptic involutions of $G$ are in 1-1 correspondence with the elliptic subfields of $K$ of degree 2 (by the Riemann-Hurwitz formula).  For the number of elliptic involutions of a genus 2 curve see \cite{deg2}. 

%If $z_1$ is an elliptic involution and $z_0$ the hyperelliptic one, then $z_2:=z_0\, z_1$ is another elliptic involution. So the elliptic involutions come naturally in pairs. This pairs also the elliptic subfields of $K$ of degree 2. Two such subfields $E_1$ and $E_2$ are paired if and only if $E_1\cap k(X)=E_2\cap k(X)$. $E_1$ and $E_2$ are $G$-conjugate unless $G\iso D_6$ or $G\iso V_4$.

%****************************************************************************************************
%%\newpage
\subsection{Automorphisms in terms of theta functions}

%*************************
%\subsection{Curves with   automorphisms }
The locus $\L_2$ of genus 2 curves $\X$ which have an elliptic involution is a closed subvariety of $\mathcal M_2$.
\begin{lem}
Let $\X$ be a genus 2 curve. Then  $Aut(\X)\iso V_4$ if and only if the theta functions of $\X$ satisfy
\begin{scriptsize}
\begin{equation}\label{V_4locus1}
\begin{split}
(\theta_1^4-\theta_2^4)(\theta_3^4-\theta_4^4)(\theta_8^4-\theta_{10}^4)
(-\theta_1^2\theta_3^2\theta_8^2\theta_2^2-\theta_1^2\theta_2^2\theta_4^2\theta_{10}^2+\theta_1^4\theta_3^2\theta_{10}^2+ \theta_3^2\theta_2^4\theta_{10}^2)\\
(\theta_3^2\theta_8^2\theta_2^2\theta_4^2-\theta_2^2\theta_4^4\theta_{10}^2+\theta_1^2\theta_3^2\theta_4^2\theta_{10}^2-\theta_3^4\theta_2^2\theta_{10}^2) (-\theta_8^4\theta_3^2\theta_2^2+\theta_8^2\theta_2^2\theta_{10}^2\theta_4^2+\theta_1^2\theta_3^2\theta_8^2\theta_{10}^2-\theta_3^2\theta_2^2\theta_{10}^4)\\
(-\theta_1^2\theta_8^4\theta_4^2-\theta_1^2\theta_{10}^4\theta_4^2+\theta_8^2\theta_2^2\theta_{10}^2\theta_4^2+\theta_1^2\theta_3^2\theta_8^2\theta_{10}^2) (-\theta_1^2\theta_8^2\theta_3^2\theta_4^2+\theta_1^2\theta_{10}^2\theta_4^4+\theta_1^2\theta_3^4\theta_{10}^2-\theta_3^2\theta_2^2\theta_{10}^2\theta_4^2)\\
(-\theta_1^2\theta_8^2\theta_2^2\theta_4^2+\theta_1^4\theta_{10}^2\theta_4^2-
\theta_1^2\theta_3^2\theta_2^2\theta_{10}^2+\theta_2^4\theta_4^2\theta_{10}^2) (-\theta_8^4\theta_2^2\theta_4^2+\theta_1^2\theta_8^2\theta_{10}^2\theta_4^2-\theta_2^2\theta_{10}^4\theta_4^2+\theta_3^2\theta_8^2\theta_2^2\theta_{10}^2)\\
(\theta_1^4\theta_8^2\theta_4^2-\theta_1^2\theta_2^2\theta_4^2\theta_{10}^2-\theta_1^2\theta_3^2\theta_8^2\theta_2^2+\theta_8^2\theta_2^4\theta_4^2) (\theta_1^4\theta_3^2\theta_8^2-\theta_1^2\theta_8^2\theta_2^2\theta_4^2-\theta_1^2\theta_3^2\theta_2^2\theta_{10}^2+\theta_3^2\theta_8^2\theta_2^4)\\
(\theta_1^2\theta_8^4\theta_3^2-\theta_1^2\theta_8^2\theta_{10}^2\theta_4^2+\theta_1^2\theta_3^2\theta_{10}^4-\theta_3^2\theta_8^2\theta_2^2\theta_{10}^2)
(\theta_1^2\theta_8^2\theta_4^4-\theta_1^2\theta_3^2\theta_4^2\theta_{10}^2+\theta_1^2\theta_3^4\theta_8^2-\theta_3^2\theta_8^2\theta_2^2\theta_4^2)
 & =0
\end{split}
\end{equation}
\end{scriptsize}
\end{lem}
This was done in \cite{theta-1}.
We would like to express the conditions of the previous lemma in terms of the fundamental theta constants only.   
\begin{lem}
Let $\X$ be a genus 2 curve. Then we have the following:
\begin{description}
\item [i] $V_4 \hookrightarrow Aut(\X)$ if and only if the fundamental theta constants of $\X$ satisfy
%%%
\begin{small}
\begin{equation}\label{V_4locus2}
\begin{split}
\left( \theta_3^4-\theta_4^4 \right)  \left(\theta_1^4 -\theta_3^4 \right) \left(
\theta_2^4-\theta_4^4 \right) \left( \theta_1^4 -\theta_4^4 \right)  \left(
\theta_3^4-\theta_2^4 \right) \left( \theta_1^4-
\theta_2^4 \right) \\
\left( -\theta_4^2+\theta_3^2+\theta_1^2-\theta_2^2 \right)\left(
 \theta_4^2-\theta_3^2+\theta_1^2-\theta_2^2
 \right)  \left( -\theta_4^2-\theta_3^2+\theta_2^2+\theta_{{
1}}^2 \right)  \left( \theta_4^2+\theta_3^2+\theta_2^2+\theta_ 1^2 \right)\\
\left( {\theta_1}^4{\theta_2}^4+
{\theta_3}^4{\theta_2}^4+{\theta_1}^4{\theta_3}^4-2\,\theta_1^2\theta_2^2\theta
_3^2\theta_4^2 \right) \left( -{\theta_3}^4{\theta_2}^4-{
\theta_2}^4{\theta_4}^4-{\theta_3}^4{\theta_4}^4 +
2\,\theta_1^2\theta_2^2\theta_
3^2\theta_4^2 \right)\\
 \left( {\theta_2}^4{\theta_4}^4+{\theta_1}^4{\theta _2}^4+{
\theta_1}^4{\theta_4}^4-2\,\theta_1^2\theta_2^2\theta_3^2\theta_4^2 \right)
\left( {\theta_1}^4{\theta_4}^4+{\theta_3}^4{\theta_4}^4+{\theta_1}^4{\theta_{
3}}^4
-2\,\theta_{{1 }}^2\theta_2^2\theta_3^2\theta_4^2\right)  = & 0\\
\end{split}
\end{equation}
\end{small}
\item [ii] $D_8 \hookrightarrow Aut(\X)$ if and only if the fundamental theta constants of $\X$ satisfy
Eq.~(3) in \cite{theta-1}

\item [iii] $D_6 \hookrightarrow Aut(\X)$ if and only if the fundamental theta constants of $\X$ satisfy
Eq.~(4) in \cite{theta-1}
\end{description}
\end{lem}

Hence, one can determine very easily the automorphism group of the curve once its thetanulls are known. However, this is not quite easy since thetanulls are defined in terms of the complex integrals.  However, since the Siegel modular  forms are expressed in terms of the thetanulls as in Eq.~\eqref{Siegel}, we should be able to express such loci in terms of the modular forms.

Instead, it is much easier to determine the automorphism groups in terms of the algebraic invariants which are given in terms of the coefficients of the curve, but are invariant of the equation of the curve.  In the next section we will express such loci in terms of the absolute invariants $i_1, i_2, i_3$. Then, using formulas in Eq.\eqref{i-modular} we can always express these equations in terms of the modular forms. 

%**********************

\subsection{Automorphisms in terms of algebraic invariants}

The locus of genus 2 curves with an extra involution is computed in \cite{deg2}*{Theorem 3} in terms of $J_2, J_4, J_6, J_{10}$.  It corresponds precisely to the locus $\det M =0$.
\begin{prop} Let $\p \in \M_2$.  The following hold:

\begin{enumerate}

\item  If $\B (\p)  \neq 0$, then $|\Aut (\p)|=2$ or $\p$ correspond to the curve $y^2=x(x^5-1)$.  Moreover,
 $\B (\p)=0$ if and only if  $V_4 \emb  \Aut (\p) $.

 \item  If $\B (\p)  = 0$ and $\dd (\p) =0$, then $V_4 \emb \bAut (\p)$.
 % is isomorphic to $D_4$ or $D_6$.  
 
 %Moreover,

% i) if  $v^2=4u^3$ then $\Aut (\p) \iso D_4$

% ii)  if ${u}^2-110\,u -4\,v+1125=0$ then $\Aut (\p) \iso D_6$.

\end{enumerate}
\end{prop}
Part i)  was known to Clebsch and proved independently in \cite{deg2}*{Theorem 3}. The proof of 2) can be found in \cite{deg2}*{Lemma 3}.
In the following proposition we give conditions on the  absolute invariants $i_1, i_2$ and $i_3$ for each automorphism to happen.

\begin{prop}
Let $\p =( i_1,  i_2,  i_3) \in \M_2 \setminus\{ J_2=0 \}$.  Then the following hold:

\begin{itemize}
\item[i)]     If $\p=(0, 0, 0)$ then $\Aut (\p) \iso C_{10}$
%        return [10,1];

\item[ii)]    If $\p=(\frac {81} {20}, - \frac {729} {200},  \frac {729} {25600000} )$, then $\Aut (\p) \iso SL_2 (3)$
%        return [24, 8];

\item[iii)]    If $\p=( - \frac {36} {5}, \frac {1512} {25}, \frac {243} {200000})$, then $\Aut (\p) \iso GL_2 (3)$
%        return [48,5];

\item[iv)]     If the following conditions are satisfied 
\begin{small} 
\begin{equation}
\begin{split}
27 i_2^4-9 i_2^4 i_1+9459597312000 i_3^2 i_1^2-111451255603200 i_3^2 i_1+55240704 i_3 i_1^4 \\
 -161243136 i_3 i_1^3+27 i_1^6+240734712102912 i_3^2+264180754022400000 i_3^3 \\
-9 i_1^7+18 i_2^2 i_1^4 -54 i_1^3 i_2^2 -161243136 i_3 i_2^2+12441600 i_3 i_2^3-107495424 i_3 i_2 i_1^2 \\
+52254720 i_3 i_2^2 i_1-2 i_2 i_1^6 +4 i_2^3 i_1^3-20639121408000 i_3^2 i_2+8294400 i_3 i_2^2 i_1^2\\
-331776 i_3 i_1^5 +2866544640000 i_3^2 i_1 i_2+47278080 i_3 i_1^3 i_2-2 i_2^5 & =0 \\
-243 i_1^2+80 i_1^3-1458 i_2+540 i_2 i_1+100 i_2^2 & =0 \\
\end{split}
\end{equation}
\end{small}
 and $\p$ is not one of the cases ii) and iii),  then $\Aut (\p) \iso D_4$
%        return [8, 3];

\item[v)]    If the following conditions are satisfied 
\begin{small} 
\begin{equation}
\begin{split}
3888 i_1+432 i_2-1188 i_1^2+5 i_1^3-25 i_2^2-360 i_2 i_1  & = 0 \\
26873856 i_3+5184000 i_3 i_1^2-9331200 i_3 i_1-149299200000 i_3^2-729 i_1^2 - \\
27 i_1^4 +243 i_1^3+i_1^5  & = 0\\
\end{split}
\end{equation}
\end{small}

 and $\p$ is not one of the cases ii) and iii),  then $\Aut (\p) \iso D_6$
%        return [12, 6];

\item[vi)]     If the following is satisfied
%i_3
\begin{small} 
\begin{equation} 
\begin{split}
-27 i_1^6-9459597312000 i_3^2 i_1^2+20639121408000 i_3^2 i_2+111451255603200 i_3^2 i_1 \\
-240734712102912 i_3^2-55240704 i_3 i_1^4-18 i_1^4 i_2^2-8294400 i_3 i_2^2 i_1^2-47278080 i_3 i_2 i_1^3\\
-264180754022400000 i_3^3-2866544640000 i_3^2 i_2 i_1+2 i_1^6 i_2-4 i_1^3 i_2^3+9 i_1^7+331776 i_3 i_1^5 \\
+107495424 i_3 i_2 i_1^2-27 i_2^4+9 i_1 i_2^4-52254720 i_3 i_2^2 i_1+2 i_2^5+161243136 i_3 i_2^2\\
+161243136 i_3 i_1^3-12441600 i_3 i_2^3+54 i_1^3 i_2^2  & =  0
\end{split}
\end{equation}
\end{small}
and $\p$ is not one of the above cases, then $\Aut (\p) \iso V_4$. 
\end{itemize}

\end{prop}

\proof The proof is computational and straightforward.  The reader can check \cite{MR2701993} or \cite{deg2} for details. 
\qed

If $\X$ has extra involutions (i.e. it has  automorphisms other then the hyperelliptic involution), then $X$
is isomorphic to a curve given by
\begin{equation}\label{eq-v4}
y^2=x^6+a x^4+b x^2 +1
\end{equation}
for appropriate values of $a,b$ (i.e. the discriminant is nonzero). Such form of genus 2 curves with automorphisms was know to XIX century mathematicians, but it has been used in the last decade quite often mostly due to Shaska and V\"olklein paper \cite{deg2} which first appeared in 2000.

In \cite{deg2} for curves with automorphism  were defined the dihedral invariants
\begin{equation}
u:  = a b,     \qquad v:  =a^3+b^3
\end{equation}
which give a birational parametrization of this locus $\L_2$ which is a 2-dimensional subvariety of $\M_2$.  It was the first time that such parametrization was discovered and since it
has become the common way of computing with genus 2 curves with automorphisms.  We can express $u$, and $v$ in terms of the absolute invariants $i_1, i_2, i_3$ as in \cite{deg2}. Moreover, the invariant $\dd$ is expressed  in terms of $u$, and $v$ as follows
\[
\dd  =     \left( 4u^3 -v^2 \right)  \left( {u}^2-110\,u -4\,v+1125 \right)^2
 \]
%
%Maybe we need to explain a bit here...............
%
We have the following:

\begin{prop}\label{u_v}
Let $\p$ be a genus 2 curve such that $G:=\Aut(\p)$ has an elliptic involution. Then,

a) $G\iso SL_2(3)$ if and only if $(u,v)=(0,0)$ or $(u,v)=(225, 6750)$.

b) $G\iso GL_2(3)$ if and only if $u=25$ and $v=-250$.

c) $G\iso D_6$ if and only if $4v-u^2+110u-1125=0$, for $u\neq 9, 70 + 30\sqrt5, 25$.

d) $G\iso D_4$ if and only if $v^2-4u^3=0$, for $u \neq 1,9, 0, 25, 225$. Cases $u=0,225$ and $u=25$ are
reduced to cases a),and b) respectively. 
\end{prop}

The $V_4$-locus (i.e., the locus $J_{30}=0$) is birationally parametrized by dihedral invariants 
 $u, v$ as in \cite{deg2}. 

%***************
\subsection{The locus $J_2=0$}   If a point $p\in \M_2$ is in the locus $J_2 =0$ then the following result holds true:

\begin{prop} Let $\p = (0, t_1, t_2, t_3) \in \M_2$. Then $| \Aut (\p) | > 2$ if and only if the invariants $t_2, t_3$ satisfy the Eq.~\eqref{eq_J2} in the Appendix. 
Moreover, $\Aut ( \p) \iso D_4$ if and only if 
\[ y^2= x^5 + x^3  - \frac 3 {20} x \]
and 
$\Aut ( \p) \iso D_6$ if and only if 
\[ y^2= x^6  + x^3 + - \frac 1 {40} \]

\end{prop}

\proof
We let $\p$ be a curve with extra automorphisms.  Then its equation can be written as 
\[ y^2 = x^6+ax^4 + bx^2+1\]
Compute $t_2, t_3$ in terms of the dihedral invariants $u, v$. Since $J_2 =0$, then $u = 0$.  Hence, $t_2 $ and $t_3 $ are expressed only in terms of $v$.  Eliminate $v$ and we get the equation in Eq.~\eqref{eq_J2}.

Assume that $\Aut ( \p) \iso D_4$.  Then, the curve has equation 
\[ y^2 = x^5+x^3+sx,  \]
cf. section 8.2. We then have 
\[ J_2 =6 + 40 s    = 0 \]
Hence,  $s= - \frac 3 {20}$.  

If $\Aut ( \p) \iso D_4$, then the equation of the curve is 
\[ y^2 = x^6 + x^3 + w, \]
cf. section 8.3. We then have 
\[ J_2 = 6 - 240 w    = 0 \]
Hence,  $ w= - \frac 1 {40}$.  

\qed

In the next section we will focus on genus two curves over  $\Q$.

%*************************************************
%%\newpage
\section{Genus 2 curves defined over \texorpdfstring{$\Q$}{Q}}

Let $\X$ be a genus two curve defined over $\Q$.  The moduli point in $\M_2$ corresponding to $\X$ is given by $\p = (i_1, i_2, i_3)$.  Since $i_1, i_2, i_3$ are rational functions in terms of the coeffiecients of $\X$, then $i_1, i_2, i_3 \in \Q$.   

The converse isn't necessarily true.  Let $\p= (i_1, i_2, i_3)  \in \M_2 (\Q)$.  The universal equation of a genus 2 curve  corresponding to $\p$ is determined  in \cite{MS-1}, which is defined over a quadratic number field $K$.  The main questions we want to consider is what percentage of the rational moduli points are defined over $\Q$?  How can we determine a minimal equation for such curves?
% How do we get an equation of such curves with minimal discriminant?  How many twists of minimal height correspond to a given $\p$?  

%****************************************************************************************************
\def\x{\mathbf{x}}
 
The quotient space $C:=\X / \Aut (\X)$ is a genus zero curve, i.e. is isomorphic to a conic which is determined as follows.  Let $\x = (x_1, x_2, x_3)$ be a coordinate in $\P^2$ and  
\begin{equation}\label{conic} C: \, \, \psi (x_1, x_2, x_3):=\x^t \, M \, \x \end{equation} 
be the corresponding conic for the symmetric matrix $M$ given as  in Eq.~\eqref{M-matrix}.    
%This is called Mestre's conic  and it is easy to see that its coefficients can be expressed in terms of the absolute invariants. 
%The conic $C$ is given by 
%
%\[ \sum_{i\leq j, j \leq k}  q_{i, k} x_i x_k =0 \]
%
%where the matrix associated to $C$ is exactly the Clebsch matrix $M$ in \eqref{M-matrix}.

There is also a cubic $L$  given by the equation 
\begin{equation}\label{cubic} L: \;  \sum_{1 \leq j, k, l\leq 3} a_{jkl} x_j x_k x_l=0\end{equation}
where the coefficients $a_{jkl}$ are given explicitly in terms of the invariants in \cite{MS-1}. 

\begin{lem}[\cite{Me}]
The genus 2 curve $\X$ is the intersection of the conic $C$ and the cubic $L$. 
\end{lem}

Mestre's algorithm has been implemented in Magma, Maple, Sage.    For an implementation in Sage see Bouyer/Streng and their paper  \cite{streng}.   However, all these implementations have had issues. For example, the implementation in Sage will not work for curves with $J_2=0$.  

\begin{exa}\label{Me-fail}
Let $\X$ be the genus 2 curve given by the equation
\[ y^2 = x^6+x^5+x+1/6. \]
One can check that for this curve the Igusa-Clebsch invariants are 
\[ [ \, 0, \, -32000,\,  5120000/3, \, 295116800000/81\, ]\]
and the algorithm fails in this case. 
\end{exa}

%Next we have our main result:
%
\begin{thm}\cite{MS-1}*{Thm. 2}
For every point  $\p \in \M_2\setminus \{ D=0 \}$ such that $\p \in \M_2 (k)$, for some number field $K$,  there is a pair of genus-two curves $\mathcal{C}^\pm$ given by
\[ \mathcal{C}^\pm: \quad  y^2= \sum_{i=0}^6 \, a_{6-i}^\pm\, x^i \; , \]
 corresponding to $\p$, such that   $a^\pm_i \in K( d )$, $i=0, \dots , 6$ as given explicitly in \cite{MS-1}*{Eq. 41}. 
% Moreover, $K ( d )$ is the minimal field of definition of $\p$. 
\end{thm}

\proof
The proof is exactly the same as \cite{MS-1}*{Thm. 4}, other then the fact that the coefficients are expressed in terms of $i_1, i_2, i_3$. We can take $J_2, \dots , J_{10}$ as in Eq.~\eqref{weighted_J}. Substituting in formulas Eq.~\eqref{iClebsch_invariants} we get, 
\begin{equation}\label{Clebsch_i}
\begin{split}
A  & =  -  \frac 1 {120}   \\
B  & =  \frac {1} {2^5\cdot 3^5 \cdot 5^4} \, (36 + 5 i_1)   \\
C &  =  \frac 1 {2^8 \cdot 3^8 \cdot 5^6}  \,   (25 i_2 +180 i_1 -216)  \\
D  & =   \frac 1 {2^{12} \cdot 3^{14} \cdot 5^{10} }  \, (2700 i_2 - 675 i_1^2 - 250 i_1 i_2 - 86400000 i_3 + 13500 i_1 - 34992 ) \\
\end{split}
\end{equation}
Then $d^2$ in \cite{MS-1} becomes 
\begin{small}
\begin{equation}\label{eq-d-2}
\begin{split}
d^2  = &  - \frac 1 {2^{50} \cdot 3^{56} \cdot 5^{30} } \, (9 i_1^7 + 2 i_1^6 i_2 - 27 i_1^6 - 18 i_1^4 i_2^2 - 4 i_1^3 i_2^3 + 331776 i_1^5 i_3 + 54 i_1^3 i_2^2 + 9 i_1 i_2^4 \\
& + 2 i_2^5 - 55240704 i_1^4 i_3 - 47278080 i_1^3 i_2 i_3 - 8294400 i_1^2 i_2^2 i_3 - 27 i_2^4 + 161243136 i_1^3 i_3 \\
& + 107495424 i_1^2 i_2 i_3 - 52254720 i_1 i_2^2 i_3 - 12441600 i_2^3 i_3 - 9459597312000 i_1^2 i_3^2 \\
& - 2866544640000 i_1 i_2 i_3^2 + 161243136 i_2^2 i_3 + 111451255603200 i_1 i_3^2 \\
& + 20639121408000 i_2 i_3^2 - 264180754022400000 i_3^3 - 240734712102912 i_3^2) \cdot \\
& \cdot (675 i_1^2 + 250 i_1 i_2 - 13500 i_1 - 2700 i_2 + 86400000 i_3 + 34992) 
\end{split}
\end{equation}
\end{small}
Notice that $d^2$ has two significant factors: one is $J_{30}$ which correspond exactly to the locus of the curves with extra involutions, and the other one is the Clebsch invariant $D$ expressed in terms of $i_1, i_2, i_3$.  
Substituting $A, B, C, D$ for the expressions in Eq.~\eqref{Clebsch_i} we get $a_0, \dots , a_6$ in terms of $i_1, i_2, i_3$ as in the Appendix. 
\qed

\begin{rem}
Notice that the above theorem does not provide an equation for the curve if the Clebsch invariant $D=0$. 
\end{rem}

%If the conic has a rational point  then we can find infinitely many rational points. Hence, if the conic has a rational point we can parametrize, say  
%
%\begin{equation}
% x_1 =   h_1(t), \quad x_2=  h_2(t), \quad x_3=  h_3(t)
% \end{equation}
%
%where $h_i \in \Q(x, y, z)[t]$ are quadratic equations in terms of the parameter $t$.  Hence, the set of points of the conic can be parametrized by the variable $t$. 

 %**********************
 
\subsection{Genus 2 curve with $\Aut( \X) \cong V_4$}

The above theorem is valid for any point $\p \in \M_2$, including the points $\p\in \M_2$ with $V_4 \embd \Aut (\p)$. On contrary, Mestre's approach which has been implemented in many computer algebra packages is valid only for $ \p \in \M_2$ with $\Aut (\p) \iso C_2$.

%However, another more direct approach for all the cases $V_4 \embd \Aut (\p)$ has been used in \cite{arith-gen-2}.

 Let $\X$ be a genus 2 curve with $\Aut( \X) \cong V_4$. The  $V_4$ locus has dimension two and is  parametrized by the parameters $u,$ and $v$ as explained above. The rational model of a genus 2 curve with automorphism group $V_4$ given in terms of this parameters $u$, and $v$ has equation given as  follows
 \begin{small}
\begin{equation}\label{eq-V_4}
\begin{split}
f(x, z)  = & \left( {v}^{2}+{u}^{2}v-2\,{u}^{3} \right) {x}^6+2\, \left( {u}^{2}+3\,v \right)  \left( {v}^{2}-4\,{u}^{3} \right) {x}^5z +\\
&+ \left( 15\,{v}^{2}-{u}^{2}v-30\,{u}^{3} \right)  \left( {v}^{2}-4\,{u}^{3} \right) {x}^{4}z^2 +4\, \left( 5\,v-{u}^{2} \right)  \left( {v}^{2}-4\,{u}^{3} \right) ^{2}{x}^{3}z^3+\\
&+ \left( {v}^{2}-4\,{u}^{3} \right) ^{2} \left( 15\,{v}^{2}-{u}^{2}v-30\,{u}^{3} \right) {x}^{2}z^4 +2\, \left( {v}^{2}-4\,{u}^{3} \right) ^{3} \left( {u}^{2}+3\,v \right) xz^5+\\
&+ \left( {v}^{2}-4\,{u}^{3} \right) ^{3} \left( {v}^{2}+{u}^{2}v-2\,{u}^{3} \right) z^6
 \end{split}
\end{equation}
\end{small}
where $u$, and $v$ are expressed in terms of the absolute invariants $(i_1, i_2, i_3)$ as in \cite{deg2}. 
   The discriminant of the form is 
\[ \D_f= 2^{36} \, \, \left( {u}^{2}+18\,u-4\,v-27 \right) ^{2} \left( 4\,{u}^{3}-{v}^{2} \right) ^{15}{u}^{30}.\]
The case $\D_f=0$ correspond exactly to the cases when $|\Aut (\p) | > 4$ which is treated below.  Notice that the factors $2$ and $u$ have exponents $\geq 30$.  

%***************************************************** 
\subsection{Genus 2 curve with $\Aut( \X) \cong D_4$}

 Let $\X$ be a genus 2 curve with $\Aut( \X) \cong D_4$. The  $D_4$ locus is a dimension one locus parametrized by the parameter $s$, as proved in \cite{ants}*{Lemma~3}{ants}. The rational model of a genus 2 curve with automorphism group $D_4$ given in terms of this parameter $s$  is as follows 
\[ y^2={x}^5+{x}^{3}+ s\,  x \]
where $s$ can be expressed in terms of the absolute invariants $(i_1, i_2, i_3)$ as follows 
\[ 
s= - \frac 3 4  \frac {345\, i_1^{2}+50\,  i_1     i_2 -1296   i_1 -90  i_2 }  {2925  i_1 ^{2}+250\  i_1   i_2 -54000  i_1 -9450  i_2 +139968 }.
%s= - \frac{9}{4} \cdot \frac A B
\]
The discriminant   is $ \D_f=2^4 \cdot \,{s}^{3}\left( 4\,s-1 \right) ^{2}$.  
%In \cite[Lemma~4]{ants} an expression of $s$ is given in terms of $i_1, i_2, i_3$. 

%*****************************************************
\subsection{Genus 2 curve with $\Aut( \X) \cong D_6$}

Now let us consider the case when  $\X$ is  a genus 2 curve with $\Aut( \X) \cong D_6$. The  $D_6$ locus is a dimension one locus parametrized by the parameter $w$ as proved  in \cite{ants}*{Lemma~3}. The rational model of a genus 2 curve with automorphism group $D_6$ given in terms of this parameter $w$  is as follows 
\[ y^2={x}^6+{x}^{3} + w  \]
where the parameter $w$ given in terms of the absolute invariants $(i_1, i_2, i_3)$ is 
\[
w= \frac 1 4  \frac {540 i_1^2+100 i_1 i_2-1728 i_1+45 i_2}     {2700 i_1^2+1000 i_1 i_2+204525 i_1+40950 i_2-708588}
\]
The discriminant of the form is   $ \D_f= - 3^6 \cdot {w}^{2} \left( 4\,w-1 \right) ^{3}$. 
%In \cite[Lemma~4]{ants} an expression of $w$ is given in terms of $i_1, i_2, i_3$. 

Hence, for every point $\p \in \M_2$ we can find an equation of a curve $C$ corresponding to $\p$.  Below we give an example to illustrate how this work in our code in Sage.
\begin{exa}\label{ex-d6}
Consider the curve
\begin{Small}
\[ 
\begin{split}
y^2= &\, \,  4294967297 \, t^6 + 77309411328 \, t^5 + 579820584969 \, t^4 + 2319282339816 \, t^3 + 5218385264643 \, t^2 \\
& + 6262062317592 \, t + 3131031158771 \\
\end{split}
\]
\end{Small}
By using the functions above we find that this curve has automorphism group $D_6$. Hence, we compute $w=2^{33}$.  
The equation of a curve over $\Q$ isomorphic (over $\C$)  to our curve is
\[ y^2= x^6+x^3+2^{33} \]
In the next section we will see if we can transform this curve over $\Q$ to another curve with smaller coefficients. 
\end{exa}

In the next section we will see how to  minimize the discriminant of a genus two curve when its equation is given in Weierstrass  form.

%*****************

\section{Minimal discriminant for Weierstrass equations}   %Rachel
Let $K$ be a field with a discrete valuation $\v$  and ring of integers $\O_K$   and $C$ an irreducible, smooth, algebraic    curve of genus $g\geq 1$ defined over $K$ and function field $K (C)$. The discriminant $\dis_{C/K}$ is an important invariant of the function field of the curve and therefore of the curve. Since the discriminant is a polynomial given in terms of the coefficients of the curve, then it is an ideal in the ring of integers $\O_K$ of $K$.  The valuation of this ideal is a positive integer.
A classical question is to find an equation of the curve such that this valuation is minimal, in other words the discriminant is minimal.

When $g=1$, so that $C$ is an elliptic curve, there is an extensive theory  of the minimal discriminant ideal $\dis_{C/K}$.  Tate \cite{ta-75}  devised an algorithm how to determine the Weierstrass equation of an elliptic curve with minimal discriminant as part of his larger project of determining Neron models for elliptic curves.  Such ideas were extended for genus 2 in \cite{liu}.   Here we mostly follow  \cite{rachel}.

For a  binary form $f(x, z)$ and a matrix
$M =  \begin{bmatrix} a, b \\ c, d \end{bmatrix}$, such that  $M \in GL_2 (k)$, we have that $f^M:= f(aX + bZ, cX + dZ)$ has discriminant
\[\D (f^M) = ( \det M )^{\frac {d(d-1)} 2} \cdot \D (f)\]
This property of the discriminant is crucial in our algorithm which is explained later.

%\begin{lem}
%i) The discriminant of a degree $d$ binary form $f(x, z)\in k[x, z]$ is and $SL_2 (k)$-invariant of degree $2d-2 $.

%ii) For any $M \in GL_2 (k)$ and any degree $d$ binary form $f$ we have that
%
%\[ \D (f^M) = \left(  \det M \right)^{d (d-1) } \, \D (f) \]
%\end{lem}
%

%***********************************************
\subsection{Discriminant of a curve}
The concept of a minimal discriminant for elliptic curves was defined by Tate and others in the 1970-s; see \cite{ta-75}.    Such definitions and results were generalized by Lockhart in \cite{lockhart} for hyperelliptic curves.

%\subsection{Minimal discriminants over local fields}
%
%Let $K$ be a local field, complete with respect to a valuation $\v$.  Let $\O_K$ be the ring of integers of $K$, in other words
%
%$\O_K = \{ x\in K \, | \, \v (x) \geq 0\}$.
%
%We denote by $\O_K^\ast$ the group of units of $\O_K$ and by $\m$ the maximal ideal of $\O_K$.  Let $\pi$ be a generator  for $\m$ and $k=\O_K / \m$ the residue field. We assume that $k$ is perfect and denote its algebraic closure by $\bar k$.

%Let $\X_g$ be a superelliptic curve of genus $g \geq 2$ defined over $K$ and $P$ a $K$-rational point on $\X_g$.     By a suitable change of coordinates we can assume that all coefficients of $\X_g$ are in $\O_K$.      Then, the discriminant $\D \in \O_K$.  In this case we say that the equation of $\X_g$ is \textbf{integral}.

%An equation for $\X_g$ is said to be a \textbf{minimal equation}  if it is integral and $\v (\D)$ is minimal among all integral equations of $\X_g$. The ideal $I=\m^{\v (\D)}$ is called the \textbf{minimal discriminant} of $\X_g$.

%*********************************************************
\newcommand\fa{\mathfrak a}

%\subsection{Minimal discriminants over global fields}
%
Let us assume now that $K$ is an algebraic number field with field of integers $\O_K$.  Let $M_K$ be the set of all inequivalent absolute values on $K$  and $M_K^0$ the set of all non-archimedean absolute values in $M_K$.
We denote by $K_\v$ the completion of $K$ for each $\v \in M_K^0$ and by $\O_\v$ the valuation ring in $K_\v$. Let $\p_v$ be the prime ideal in $\O_K$ and $\m_v$ the corresponding maximal ideal in $K_\v$. Let $(\X, P)$ be a superelliptic curve of genus $g\geq 2$ over $K$.

If $\v \in M_K^0$ we say that $\X$ is \textbf{integral at $\v$} if $\X$ is integral when viewed as a curve over $K_\v$.  We say that $\X$ is \textbf{minimal at $\v$} when it is minimal over $K_\v$.

An equation of $\X$ over $K$ is called \textbf{integral} (resp. \textbf{minimal}) over $K$ if it is integral (resp. minimal) over $K_\v$, for each $\v \in M_K^0$.

Next we will define the minimal discriminant over $K$ to be the product of all the local minimal discriminants. For each $\v \in M_K^0$ we denote by $\D_\v$ the minimal discriminant for $(\X, P)$ over $K_\v$.  The \textbf{minimal discriminant} of $(\X, P)$ over $K$ is the ideal
\[ \D_{\X / K} = \prod_{\v \in M_K^0} \m_\v^{\v (\D_\v) } \]
We denote by $\fa_\X$ the ideal   $ \fa_\X = \prod_{\v \in M_K^0} \p_\v^{\v (\D_\v) }$.

%*****************************************************************************
%\subsection{Superelliptic curves with minimal discriminant}
%
Let $\X$ be a genus two    curve with equation $y^2=f(x)$. 
The discriminant of $\X$ is the discriminant of  $f(x)$, hence $J_{10} (f)$, a degree 10 polynomial in terms of the coefficients of $f(x)$. 

Let $M \in GL_2 (K)$ such that $\det M = \l$. Then,   $ \D (f^M) =   \l^{30} \,  \D(f)$.
Assume $\D (f) = p^\alpha \cdot N$, where $\a>30$,  for some prime $p$ and some integer $N$ such that $(p, N)=1$. 
We perform the coordinate change 
\[  x \to \frac 1 {p} x\]
 on  $f(x)$.  Then,  the new discriminant is  $\D^\prime = \frac 1 {p^{30} } \, \cdot  p^\alpha \cdot N = p^{\alpha - 30} \cdot N$.
Hence, we have the following
\begin{lem}  A genus 2 curve $\X_g$ with integral equation
\[ y^2 = a_6 x^6 + \cdots a_1 x + a_0 \]
has minimal discriminant  if $\v (\D) <  30$.
\end{lem}
Thus, we factor $\D$ as a product of primes, say $\D = p_1^{\alpha_1} \cdots p_r^{\alpha_r}$,   and take $u$ to be the product of those powers of primes with exponents $\alpha_i \geq 30$.  
Then, the transformation $x \mapsto \frac 1 {p^u} \cdot x$ will reduce the discriminant.

%***********************************
\subsection{Twists }
Of course we can make the discriminant as small as we want it if we allow twists.

Let $\X$ be hyperelliptic curve with an extra automorphism of order $n\geq 2$.  Then, from \cite{sh-issac}  we know that the equation of $\X$ can be written as 
\[ y^2= f(x^n) \]
If $\X$ is given with such equation over $\Q$ and discriminant 
\[  \D = u \, N,  \quad \textit{such that } \quad (u, N)=1  \]
then for any transformation $\tau:  x \mapsto  {u^{- \frac n {30}}}  \cdot x$ 
would lead to  $\X^\tau$ defined over $\Q$ and isomorphic to $\X$ over $\C$. Hence, $\X^\tau$ is a twist of $\X$ with discriminant 
\[ \D^\prime =  \frac 1 {u } \cdot  u \, N = N, \]
Hence, in this case we can reduce the discriminant even further.

\begin{lem}
Let $\X$ be a genus 2 curve with $\Aut (\X) \iso V_4$.  Assume that $\X$ has equation over $\Q$ as $y^2= f(x^2)$.  Then, $\X$ has minimal discriminant $\D$ if $\v (\D) <  15$.

If $\Aut (\X) \iso D_6$ and    $\X$ has equation over $\Q$ as $y^2= f(x^3)$.  Then, $\X$ has minimal discriminant $\D$ if $\v (\D) <  10$.
\end{lem}

\proof   The proof is an immediate consequence of the above remarks. \qed

The following gives a universal curve with minimal discriminant for the $V_4$-locus.
\begin{lem}
Let $\X$ be a genus 2 curve such that $\Aut (\X)\iso V_4$ and $(u, v)$ the corresponding dihedral invariants.  Then, $\X$ has minimal discriminant and is defined over its field of moduli when given by the equation
\begin{equation}\label{min_disc}
 y^2=b_6x^6+ \cdots + b_1 x + b_0  
\end{equation} 
where 
\[
\begin{split}
b_6 = &  - {\frac { \left( 2\,{u}^{3}-{u}^{2}v-{v}^{2} \right) }{2^6 {u}^{6}}}  \\
b_5 = & - 2 {\frac { \left( {u}^{2}+3\,v \right)  \left(  4\,{u}^{3}-{v}^{2} \right) }{ 2^5 {u}^{5}}}  \\
b_4 = & {\frac { \left( 30 \,{u}^{3}+{u}^{2}v-15\,{v}^{2} \right)  \left( 4\,{u}^{3}-{v}^{2} \right) }{2^4 {u}^{4}}}  \\
b_3 = & - 4\,{\frac { \left( {u}^{2}-5\,v \right)  \left( 4\,{u}^{3}-{v}^{2} \right) ^{2}}{2^3 {u}^{3}}}  \\
b_2 = & - {\frac  { \left( 30\,{u}^{3}+{u}^{2}v-15\,{v}^{2} \right)  \left( 4\,{u}^{3}-{v}^{2} \right) ^{2}}{2^2 {u}^{2}}}   \\
b_1 = & -{\frac { \left( {u}^{2}+3\,v \right)  \left( 4\,{u}^{3}-{v}^{2} \right) ^{3}}{u}}   \\
b_0 = & \left( 2\,{u}^{3}-{u}^{2}v-{v}^{2} \right)  \left( 4\,{u}^{3}-{v}^{2} \right)^{3}  \\
\end{split}
\]
\end{lem}

\proof 
We know that the equation of the curve in Eq.~\eqref{eq-V_4} is defined over its field of moduli.  Let 
\[ \sigma:  (x, y) \mapsto \left( \frac 1 {2u} \, x, \frac 1 {(2u)^6} \, y\right).\]
Then the discriminant of the new curve $\X^\sigma$ is 
\[ \D_f= 2^{6} \, \, \left( {u}^{2}+18\,u-4\,v-27 \right) ^{2} \, \, \left( 4u^3-v^2 \right)^{15}.\]
Since all exponents are $<  15$, this discriminant can not be further reduced over $\Q$. The equation of the curve in this case becomes as in Eq.~\eqref{min_disc}. 

\qed 

Notice that for any curve $\X$ with  $\Aut (\X)\iso V_4$, we have $u \neq 0$, so the curve in  Eq.~\eqref{min_disc} is defined everywhere.  

In the next section we will see that smaller discriminant doesn't necessarily mean small coefficients.  We illustrate below with an example.

\begin{exa}  Let $u=35$ and $v=6$. Using the function 
\begin{verbatim}RatModSha_uv (35, 6)\end{verbatim}
we get 
\begin{Small}
\begin{verbatim}
(-4)*(19591*t^6+106564876*t^5-55428309960*t^4+35132884438720*t^3
+9503959738981440*t^2+3132997049150231296*t-98758721142240034304)
\end{verbatim}
\end{Small}

Notice that this curve has large coefficients. Especially, when we have that for $u=35$ and $v=6$ we get corresponding $(a, b)=(2, 3)$ and therefore a curve 
\[ y^2 = x^6+2x^4+3x^3+1 \]
As it is shown in \cite{MR3525576}, it is almost always true that the curve $y^2= x^6+ax^4+bx^2+1$ has smaller coefficients when such  $a, b \in \Q$ do exist. 
\end{exa}

%********************

%\def\L{\mathcal L}

\section{Constructing the databases}\label{data}

Now that we have all the necessary results we are ready to construct all three dictionaries.

\subsection{Curves with height $\leq 10$: the dictionary $\L_1$}
In order to create the list $\L_1$ we first compile a list of all 7-tuples $(a_0, \dots, a_6)$ and from this list eliminate all the tuples with $J_{10}=0$.  Then, for each tuple we compute the moduli point $\p = (i_1, i_2, i_3)$ and all the other invariants as described in the previous sections.

\begin{table}[hb]
\caption{The number of curves for a given height}
\begin{center}
\begin{tabular}{|c|c|c|c|c|c| }
\hline
$\h$ & $\P^6$    &   $\#$ of curves & $V_4$ & $D_4$ & $D_6$   \\
\hline
1    &  1 093     &  230     & 28    & 11  & 2   \\
% \hline
2    &  37 969    &  8 593    & 230   & 40  & 7    \\
% \hline
3    &  409 585  & 88 836     & 1054  & 112 & 26 \\
%4    &  2 351 359 & 472 040    &       &     &    \\
 \hline
\end{tabular}
\end{center}
\label{tab-1}
\end{table}%

\begin{rem}
The bound for the number of  genus 2 curves of height $\h$ is $< (2h+1)^7$.  For example, for $\h=1$ this bound  is $< 3^7 = 2187$.  In fact, as shown in Table~\ref{tab-1} the number of such curves is 230, all of which are listed in \cite{MR3525576}*{Table 2, 3}.
\end{rem}

In the Table~\ref{tab-1} we display  the  number of curves with extra automorphism for a given height $\h$.

%***************************
\subsection{Curve with extra involutions: the dictionary $\L_2$}
The list $\L_2$ was computed in \cite{MR3525576} for $1 \leq \h \leq 100$. Since every genus 2 curve with extra automorphisms can be written (over $\C)$ as 
\begin{equation}\label{v4}
y^2 = x^6 + ax^4 + bx^2 +1
\end{equation}
we create the list of tuples $(1,0, b, 0, a, 0, 1)$ for $a, b \leq \h$.  We go through the lists and delete the ones which have $J_{10}=0$. This number is displayed in  
\cite{MR3525576}*{Table 1} for $1 \leq \mh \leq 100$, including information on how many points from $\L_2$ are already in $\L_1$, how many of these curves have automorphism group $D_4$ and how many have automorphism group $D_6$.

 \begin{table}[ht]
\caption{Curves with extra involutions and height $\h \leq 101$}
\begin{center}
\begin{tabular}{|c|c|c|c|c|c|}
\hline
$\h$ &                   &  new       pts.                & $D_4$   & $D_6$ & Total pts       \\
     &   $J_{10}\neq 0$  &  in $\M_2$                     &       &      &     \\

\hline
1 & 8 & 4 & 1 & 0 & 5  \\
2 & 24 & 9 & 3 & 0 & 14  \\
3 & 47 & 12 & 4 & 0 & 26  \\
4 & 79 & 17 & 6 & 0 & 43  \\
5 & 119 & 20 & 7 & 0 & 63  \\
6 & 167 & 25 & 9 & 0 & 88  \\
7 & 223 & 28 & 11 & 0 & 116  \\
8 & 287 & 33 & 13 & 0 & 149  \\
9 & 359 & 36 & 15 & 0 & 185  \\
10 & 439 & 41 & 17 & 0 & 226  \\
\hline 
\end{tabular}
\end{center}
\label{tab-2}
\end{table}%

In \cite{MR3525576} it is discussed when such curves have minimal height and how a reduction as in \cite{MR3525574} is easier to perform in this case. In Table~\ref{tab-2} we display only the first ten rows of Table~1 from \cite{MR3525576}.

There are 20 697 curves in $\L_2$, such that for each $\h$ we have roughly $4\h$ curves. So it is expected that the number of curves of height $\leq h$,  with equation Eq.~\eqref{v4}, defined over $\Q$ is $\leq 4 \frac {\h (\h+1)} 2$; see Beshaj \cite{MR3525576} for more details.  
% Notice that such curves when defined over $\Q$ have very nice minimal equations as pointed out in \cite{MR3525576}

The first curve with automorphism group isomorphic to $D_4$ occurs for $\h = 1$. It is only one such curve,  namely the curve with equation 
\[ y^2= x^6 +x^4+x^2 +1. \]
The are only two curves with automorphism group isomorphic to $D_6$, which occur for $\h = 79$ and $\h=83$.  namely the curve  
\[ y^2=x^6+79x^4-17x^2+1 \]
and the curve
\[ y^2 = x^6+83x^6+19x^2+1.\]
There are 195 curves with automorphism group isomorphic to $D_4$.  For the largest height $\h = 101$ there are 405 curves two of  which have  automorphism group isomorphic to $D_4$.

In \cite{MR3525574} is proved the following theorem:
\begin{thm}[\cite{MR3525574}]
Let ${\mathfrak p} \in \mathcal M_2 (\Q)$ be such that $\Aut({\mathfrak p}) \iso V_4$. There is a genus 2 curve  $\mathcal X$ corresponding to $\mathfrak p$ with equation 
$y^2\, z^4 =f(x^2, z^2)$,  where 
\begin{equation}\label{eq_form}
 f(x, z)= x^6 - s_1 x^4 z^2 + s_2 x^2 z^4 - z^6. 
\end{equation} 
If $f \in \Z [x, z]$, then $f(x, z)$ or $f(-z, x)$ is a reduced binary form. 
%and $\Delta_{f(x, z) }$ is minimal over $\Z$, then $\mathcal X$ is reduced in its orbit. 
%Then there exists a representative curve $\mathcal X: y^r z^{n-r}= f(x, z)$ of ${\mathfrak p}$  such that $f$ is a totally real form. Moreover, 
\end{thm}
It is interesting to not that from 20 292 such curves we found only 57 which do not have minimal absolute height. For more details see the discussion in the last section of \cite{MR3525574}.

%********************************
\subsection{Curves with small moduli height: the dictionary $\L_3$}
For the last dictionary $\L_3$ we create a list of all points $[x_0: x_1: x_2: x_3]$  of projective height $\leq \mh$ in $\P^3 (\Q)$, for some integer $\mh \geq 1$.  Each such point correspond to the point $[J_2^5 : J_4 J_2^3 : J_6 J_2^2 : J_{10}]$.  The number of such tuples (up to equivalence in $\P^3$)   is given in column 2 of Table~\ref{tab-mod}.  Not all such points correspond to a genus 2 curve.  We delete all the points such that $x_3=0$ since they correspond to the cases when $J_{10}=0$.  What is left on the list is the number of points in the moduli space $\M_2$ with moduli height $\leq \mh$. This number is given in the third column in the Table~\ref{tab-mod}.

For each given point $\p = (r, i_1, i_2, i_3)$ in this list we find the equation of the curves  as described above.  As we already know we don't get a curve defined over $\Q$ in each case.  The number of points which are defined over $\Q$ is given in column four. From those points $\p \in \M_2$ such that the field of definition is $\Q$, the number of points with automorphisms is given in column five.

\begin{table}[h]
\caption{Genus 2 curves with bounded moduli height} 
\label{tab-mod}
\begin{center}
\begin{tabular}{|c|c|c|c|c|c|}
\hline
$\mh$   &  $\P^3 (\Q) $  & $\p \in \M_2 (\Q)$ &  $M_\p = F_\p$ &  $| \Aut (\p) | >2$  & ratio \\
        &     $n_1$           &      $n_2$              &    $n_3$            &       $n_4$               &        \\
\hline
1 & 40  & 27  &  20  &  15   &   0.25 \\
2   &  272  &   223  &   124  &   75   &   0.4 \\
3   &  1120  &   975  &   514  &   243   &   0.53  \\
4   &  2928  &   2639  &   1311  &   507   &   0.61  \\
5   &  6928  &   6351  &   3056  &   1035   &   0.66  \\
6   &  12768  &   11903  &   5561  &   1587   &   0.71  \\
7   &  23760  &   22319  &   9963  &   2667   &   0.73  \\
8   &  38128  &   36111  &   15648  &   3771   &   0.76  \\
9   &  60640  &   57759  &   24214  &   5427   &   0.78  \\
10   &  88448  &   84703  &   34936  &   7107   &   0.8 \\
11   &  131088  &   125903  &   50630  &   9867   &   0.81  \\
12   &  177712  &   171375  &   68046  &   12123   &   0.82   \\
13   &  248080  &   239727  &   93229  &   16011   &   0.83   \\
14   &  324736  &   314655  &   120358  &   19395   &   0.84   \\
15   &  427968  &   415583  &   157569  &   23907   &   0.85   \\
16   &  542720  &   528031  &   199136  &   28419   &   0.86   \\
17   &  700032  &   681887  &   252947  &   35139   &   0.86   \\
18   &  857328  &   836591  &   307964  &   40251   &   0.87   \\
19   &  1076928  &   1051871  &   381219  &   48675   &   0.87   \\
\hline
\hline
\end{tabular}
\end{center}
\end{table}%

The following question is natural: What percentage of rational points $\p \in \M_2 (\Q)$  with a fixed moduli height $\mh$ have $\Q$ as a field of definition, when $\mh $ becomes arbitrarily large?   In other words, what is the limit
\[ \lim_{\mh \to \infty } \frac {n_3} {n_2} ? \]

\begin{lem}
For large enough moduli height $\mh \in \M_2(\Q)$, the majority of genus 2 curves with moduli height $\mh$ are not defined over $\Q$.  
\end{lem}

\proof
Let $\p \in \M_2 (\Q)$ with moduli height $\mh_0$. Then, $\p = (i_1, i_2, i_3)$ for $i_1, i_2, i_3 \in \Q$.
 The equation of a curve $\X$ corresponding to $\p$ is the intersection of the conic $C : \; \psi (x_1, x_2, x_3) =0$ and the cubic $L: \; \phi(x_1, x_2, x_3)=0$ which are both defined over $\Q$, since their equations are given in terms of $i_1, i_2, i_3$.   The intersection $C \cap L$ is obtained as the resultant of the corresponding equations with respect to one of the variables $x_1, x_2, x_3$.     This resultant contains a square root which is given in terms of $i_1, i_2, i_3$; see \cite{MS-1}.   $\X$ is defined over $\Q$ if and only if the expression inside the square root is a complete square.  Since this is a surface in $\Q^3$, the set of points of fixed projective height which make it a complete square has measure  zero.  This completes the proof.  
\qed

Next we will consider a similar question.   What percentage of curves defined over $\Q$ don't have extra automorphisms when $\mh$ becomes arbitrary large?   In other words, what it the limit
\[ \lim_{\mh \to \infty } \frac {n_3 - n_4} {n_3} ? \]
Thus, we want to    determine the ratio of points of height $\mh$ in $\M_2 ( \Q)$ for which $M_{\p}\neq F_{\p}$ over the total number of points of height $\mh$ in $\M_2( \Q)$ as $\mh \to \infty$.

\begin{lem}
For large enough moduli height $\mh$, the majority of genus 2 curves with moduli height $\mh$ don't have extra automorphisms.  
\end{lem}

\proof
The number of points with height $\mh$ in the projective space increases as $\mh$ increases.  However,  such points $\p = (i_1, i_2, i_3)$ have to satisfy the equation $J_{30}=0$ when they have automorphisms.  
\qed

There is another database of genus 2 curves described in \cite{Booker:2016aa} which is has all the curves with discriminant $< 1000$.  We checked our database to see how many of our curves would be with small discriminant. We have to warn the reader that our definition of the discriminant is just $J_{10}$ and slightly different from what is used in \cite{Booker:2016aa}.  The majority of curves in the third dictionary have small discriminant, but this is not a surprise since small moduli height forces $J_{10}$ to be small.  Only one curve in the second dictionary has $ | J_{10} | < 1000$. Six curves in the first dictionary have $ | J_{10} | < 10000$.  It would be interesting to see exactly how the database in \cite{data} and \cite{Booker:2016aa} intersect.

%&&&&&&&&&&&&&&&&&&&&&&&&&&&&&&&&&&&&&&&&&&&&&
\clearpage

\appendix

\pagestyle{empty}

%************************************
%
%\pagestyle{empty}
%\newpage

%\begin{landscape}
\section{Functions of the genus 2 package}\label{app-B}

\raggedright

{\footnotesize
%\begin{multicols}{1}

% multicol parameters
% These lengths are set only within the two main columns
%\setlength{\columnseprule}{0.25pt}
%\setlength{\premulticols}{1pt}
%\setlength{\postmulticols}{1pt}
%\setlength{\multicolsep}{1pt}
%\setlength{\columnsep}{2pt}

\section*{Genus 2 package}

\noindent \textbf{Input: a sextic polynomial $f(t)$} \\

\begin{tabular}{@{}ll@{}}
\verb!J2!    &   Igusa invariant $J_2$ \\
\verb!J4!    &  Igusa invariant $J_4$ \\
\verb!J6!    &  Igusa invariant $j_6$ \\
\verb!J10!    &  Igusa invariant $J_{10}$ \\
\verb!J30!    &  $J_{30}$: the $V_4$-locus \\
\verb!Igusa!  & Igusa invariants $[J_2, J_4, J_6, J_{10}]$ \\
\verb!Clebsch!    &  Invariants [A, B, C, D] \\
\verb!i_1!    &   absolute invariant $i_1$\\
\verb!i_2!    &   absolute invariant $i_2$\\
\verb!i_3!    &   absolute invariant $i_3$\\
\verb!j1!    &   Igusa function $j_1$\\
\verb!j2!    &   Igusa function $j_2$\\
\verb!j3!    &   Igusa function $j_3$\\

\verb!RatMod!    &  Rational model of the curve over when such model exists.\\
\verb!RatModSha!    & Rational model  of the curve over its minimal field of definition as in Shaska \cite{deg2}  \\
\verb!RatModMe! &  Rational model over $\Q$, when such model exists,  as in  Mestre \cite{Me}   \\
\verb!height! &   Height of the sextic  \\
\verb!EquivBin! & Checks if sextics are equivalent    \\
\verb!RatModTable! & Rational Model from the          Table of minimal models    \\
\verb!MinField! &  Minimal field of definition   \\
\verb!Info! & Displays information about the curve  $y^2=f(t)$  \\
\verb!RatForm! &  Rational Model from     Malmendier/Shaska  \cite{MS-1} \\ 
\end{tabular}

\medskip

\noindent \textbf{Input: the moduli point $(i_1, i_2, i_3)$} \\

\begin{tabular}{@{}ll@{}}
\verb!J30_j!    & $J_{30}$ in terms of $i_1, i_2, i_3$ \\
\verb!Igusa_i!  &  $J_2, \dots , J_{10}$ \\
\verb!Clebsch_i!    &  Invariants [A, B, C, D] \\
\verb!ClebschMatrix!   & ClebschMatrix in terms of invariants $i_1, i_2, i_3$ \\
\verb!Cubic!       & Cubic defined in  Eq.~\eqref{cubic} \\
\verb!Conic!       &  Conic defined in Eq.~\eqref{conic} \\
\verb!L_D4!    & Locus of curves with group $D_4$  \\
\verb!L_D6!    & Locus of curves with group $D_6$  \\
\verb!AutGroup!  & Automorphism group of the curve    \\
\verb!Sh_u_v!    & Shaska invariants $u, v$  \\

\verb!ModHeight! & Modular height    \\ 
\end{tabular}

\section*{Moduli Space}

\begin{tabular}{@{}ll@{}}
\verb!curves_moduli!    &  Computes the number of rational points of height $\mh$ in the moduli space and\\
& how many of those have a rational model  \\
\verb!NumbCurvMod!  & This function computes the total number of rational points of moduli height $\mh$, \\
&    how many of them have a rational model over $\Q$,  how many of them have\\
& automorphisms and the ratio=obstruction/ fine points\\ 
\verb!moduli_points! &  Computes the number of rational points of height $\mh$ in the moduli space  \\
\verb!MoPtsCurvAut!  & Moduli points with automorphisms  \\        %moduli_points_of_curves_with_automorphisms
\verb!!  &  
\end{tabular}

\section*{Creating the databases}

\begin{tabular}{@{}ll@{}}
\verb!Curves(h, L)!    &  Creates the dictionary $\L_1$  of curves with height $h$  \\
\verb!CurvesAut(h, L)!  &  Creates the dictionary $\L_2$  of curves with automorphisms  \\
\verb!CurvHe! &  Number of curves with height $h$ \\
\verb!SelfRec(f)!  &  Checks if a sextic is self reciprocal \\
\verb!CurvHeW(h, w)!  & Number of curves with height $h$ and $w$ \\
\verb!NCWT(h, w)!  &  Number of curves with height $\h$ and twists $w$\\
\verb!CurvesTabOverQ(h, w)!  & Counts the number of curves over $\Q$, including  twist, for given height. \\
%\verb!curves_height(h)!  &  \\
\end{tabular}

%\end{multicols}
%\endfootnotesize
}

%***********************
\section{Basic Invariants and relations among them}\label{app-A}

\noindent Clebsch invariants $A, B, C $ are in terms of the coefficients
%
%Basic Invariants and relations among them 
\begingroup\makeatletter\def\f@size{8}\check@mathfonts\def\maketag@@@#1{\hbox{\m@th\large\normalfont#1}}%
\[
\begin{split}
 A  = & \;  2\, a_6\, a_0- \frac 1 3 \, a_5\, a_1 + \frac 2 {15}  \, a_4\, a_2- \frac 1 {20} \, a_3^2\\
 B  = & - \left(\frac 4 {225} \, a_6\, a_2\, a_4\, a_0- \frac 2 {75} \, a_6\, a_2\, a_3\, a_1+ \frac 8 {1125} \, a_6\, a_2^3- \frac 2 {75} \, a_5\, a_3\, a_4\, a_0 \frac 2 {225} \, a_5 a_3^2\, a_1 \right. \\
 & -\frac 2 {1125} \, a_5\, a_3\, a_2^2
 + \frac 8 {1125}\, a_4^3\, a_0- \frac 2 {1125}\, a_4^2\, a_3\, a_1+  \frac {14}{5625} \, a_4^2\, a_2^2- \frac 2 {225}\, a_5\, a_2\, a_4\, a_1\\
 &+ \frac2{45}\, a_5^2\, a_2\, a_0   +\frac 2 {45}\, a_6\, a_1^2\, a_4- \frac 2 9 \, a_6\, a_1\, a_5\, a_0- \frac 8 {5625} \, a_4\, a_3^2\, a_2+ \frac 2 3 \, a_6^2\, a_0^2\\
  &\left.+ \frac 2 {75}\, a_6\, a_0\, a_3^2+ \frac 1{3750}\, a_3^4  \right) \\
 C  = &  -\frac 2 {46875}\, a_4^3\, a_2^3- \frac 1 {5625}\, a_4^4\, a_1^2- \frac 1 {5625}\, a_5^2\, a_2^4-\frac 2 9 \, a_6^3\, a_0^3
 + \frac{11}{5625}\, a_6\, a_2\, a_3^2\, a_4\, a_0\\
 & - \frac7{1125}\, a_6\, a_2^2\, a_5\, a_0\, a_3+ \frac4{225}\, a_6\, a_2\, a_4\, a_1\, a_5\, a_0- \frac1{562500}\, a_3^6+\frac 8{1125}\, a_6^2\, a_2^3\, a_0\\
 &- \frac 8 {28125}\, a_6\, a_2^4\, a_4+\frac 1 {9375}\, a_6\, a_2^3\, a_3^2+ \frac 1 {5625}\, a_5\, a_3^4\, a_1+ \frac 1 {90}\, a_5^3\, a_3\, a_0^2- \frac2 {28125}\, a_5\, a_3^3\, a_2^2\\
 &- \frac 8 {28125} \, a_4^4\, a_2\, a_0 + \frac 8 {1125} \, a_4^3\, a_6\, a_0^2+ \frac 1 {9375}\, a_4^3\, a_3^2\, a_0- \frac2 {28125}\, a_4^2\, a_3^3\, a_1- \frac 1 {225}\, a_4^2\, a_5^2\, a_0^2\\
 &- \frac 2 {140625}\, a_4^2\, a_3^2\, a_2^2  + \frac 1 {90}\, a_6^2\, a_1^3\, a_3- \frac 1 {225}\, a_6^2\, a_1^2\, a_2^2+ \frac 2 {140625}\, a_4\, a_3^4\, a_2- \frac 1 {3750}\, a_6\, a_0\, a_3^4\\
 &- \frac 1 {75}\, a_6^2\, a_0^2\, a_3^2 
  +\frac {14}{225}\, a_6^2\, a_2\, a_0^2\, a_4- \frac 1 {3750}\, a_6\, a_2\, a_3^3\, a_1-\frac 2{45}\, a_6\, a_2\, a_5^2\, a_0^2- \frac1{3750}\, a_5\, a_3^3\, a_4\, a_0\\
  &+\frac 1{2250}\, a_5^2\, a_3^2\, a_0\, a_2+ \frac 7{28125}\, a_5\, a_3\, a_4\, a_2^3+\frac 1{2250}\, a_5\, a_3\, a_4^2\, a_1^2+ \frac 7{28125}\, a_4^3\, a_2\, a_3\, a_1\\
  &+\frac 2{1125}\, a_4^3\, a_1\, a_5\, a_0- \frac 1{5625}\, a_5\, a_2^2\, a_4^2\, a_1+\frac 2 {1125}\, a_5\, a_2^3\, a_6\, a_1+ \frac 1{2250}\, a_5^2\, a_2^2\, a_3\, a_1\\
  &+\frac1 9 \, a_6^2\, a_1\, a_0^2\, a_5+ \frac 1{2250}\, a_6\, a_1^2\, a_3^2\, a_4- \frac2{45}\, a_6^2\, a_1^2\, a_0\, a_4 - \frac 2 {1875}\, a_6\, a_2^2\, a_4^2\, a_0\\
  &+ \frac 1 {1875}\, a_6\, a_2^2\, a_4\, a_3\, a_1-\frac 1 {75}\, a_6^2\, a_2\, a_0\, a_3\, a_1+ \frac 1 {1875}\, a_5\, a_3\, a_4^2\, a_2\, a_0- \frac 1{75}\, a_5\, a_3\, a_6\, a_0^2\, a_4\\
  &- \frac 7 {11250}\, a_5\, a_3^2\, a_4\, a_2\, a_1+ \frac 1{90}\, a_5\, a_3^2\, a_6\, a_0\, a_1- \frac 1 {225}\, a_5^2\, a_3\, a_4\, a_1\, a_0\\
 & - \frac 7{1125}\, a_4^2\, a_6\, a_0\, a_3\, a_1-\frac 1 {225}\, a_5\, a_2\, a_6\, a_1^2\, a_3 
% D =  &   \\
\end{split}
\]

\noindent Clebsch invariants $A, B, C, D $  in terms of the Igusa invariants $J_2, J_4, J_6, J_{10}$ 
\begin{equation*}\label{iClebsch_invariants}
\begin{split}
 A  = &\; -\frac{1}{2^3 3 \cdot 5} \, J_2\;, \\
 B  = & \; \phantom{-} \frac{1}{2^3 3^3 5^4} \left( J_2^2 + 20 \, J_4 \right)\;,\\
 C = & \; - \frac{1}{2^5 3^5 5^6} \left( J_2^3 + 80 \, J_2 J_4 - 600 \, J_6\right)\;,\\
 D = & \; - \frac{1}{2^8 3^9 5^{10}} \left( 9 \, J_2^5+700 \, J_2^3 J_4-3600 \, J_2^2 J_6 -12400\, J_2 J_4^2+48000 \, J_4 \, J_6+10800000 \, J_{10}\right) \\
 \end{split}
\end{equation*} 
\noindent  The entries of the Clebsch matrix $M$.
\begin{equation*}\label{Clebsch-matrix-entries}
\begin{split}
 A_{11}  & =   2 \, C + \frac{1}{3} \, A\, B\;, \\
 A_{22}  & =  A_{13} = \, D\;,\\
 A_{33}  & =  \frac{1}{2} \, B \, D + \frac{2}{9} \, C \, (B^2+ A \, C) \;,\\
 A_{23}  & =  \frac{1}{3} \, B \, (B^2+ A \, C) + \frac{1}{3} \, C \, (2\, C + \frac{1}{3} \, A \, B) \;,\\
 A_{12}  & =  \frac{2}{3} \, (B^2+ A\, C) \\
 \end{split}
\end{equation*} 
\endgroup

Igusa invariants $J_2, J_4, J_6, J_{10}$ we display them below in terms of the coefficients

\begingroup\makeatletter\def\f@size{8}\check@mathfonts\def\maketag@@@#1{\hbox{\m@th\large\normalfont#1}}%

\begin{equation*}\label{eq_J}
 \begin{split}
J_2 = & \; 6a_3^2  -240a_0a_6+40a_1a_5-16a_2a_4\\
J_4 = &\;   48a_0a_4^3+48a_2^3a_6+4a_2^2a_4^2+1620a_0^2a_6^2+36a_1a_3^2a_5-12a_1a_3a_4^2-12a_2^2a_3 a_5+300a_1^2a_4a_6 \\
& +300a_0a_5^2a_2 +324a_0a_6a_3^2-504a_0a_4a_2a_6-180a_0a_4a_3a_5  -180a_1a_3a_2a_6+4a_1a_4a_2a_5 \\
& -540a_0a_5a_1a_6-80a_1^2a_5^2\\
J_6 = & \;   176a_1^2a_5^2a_3^2 + 64a_1^2a_5^2a_4a_2 + 1600a_1^3a_5a_4a_6 + 1600a_1a_5^3a_0a_2  - 160a_0a_4^4a_2 - 96a_0^2a_4^3a_6 \\
&+ 60a_0a_4^3a_3^2 + 72a_1a_3^4a_5  - 24a_1a_3^3a_4^2 - 119880a_0^3a_6^3  - 320a_1^3a_5^3 - 160a_2^4a_4a_6 - 96a_2^3a_0a_6^2 \\
&+ 60a_2^3a_3^2a_6 - 24a_2^2a_3^3a_5  + 8a_2^2a_3^2a_4^2  - 900a_2^2a_1^2a_6^2  - 24a_2^3a_4^3 - 36a_2^4a_5^2 - 36a_1^2a_4^4 \\
&+ 424a_0a_4^2a_2^2a_6 - 2240a_1^2a_5^2a_0a_6 + 2250a_1^3a_3a_6^2  + 492a_0a_4^2a_2a_3a_5 + 20664a_0^2a_4a_6^2a_2 \\
&+ 3060a_0^2a_4a_6a_3a_5  - 468a_0a_4a_3^2a_2a_6  - 198a_0a_4a_3^3a_5 - 640a_0a_4a_2^2a_5^2 + 3472a_0a_4a_2a_5a_1a_6 \\
&- 18600a_0a_4a_1^2a_6^2 - 876a_0a_4^2a_1a_6a_3  + 492a_1a_3a_2^2a_4a_6  - 238a_1a_3^2a_2a_4a_5 + 76a_1a_3a_2a_4^3   \\
&+ 3060a_1a_3a_0a_6^2a_2 + 1818a_1a_3^2a_0a_6a_5  + 26a_1a_3a_2^2a_5^2   - 1860a_1^2a_3a_2a_5a_6 + 330a_1^2a_3^2a_6a_4 \\
&+ 76a_2^3a_4a_3a_5  - 876a_2^2a_0a_6a_3a_5   + 616a_2^3a_5a_1a_6 + 2250a_0^2a_5^3a_3  - 10044a_0^2a_6^2a_3^2 \\
&+ 28a_1a_4^2a_2^2a_5 - 640a_1^2a_4^2a_2a_6 + 26a_1^2a_4^2a_3a_5 - 1860a_1a_4a_0a_5^2a_3   + 616a_1a_4^3a_0a_5  \\
&- 18600a_0^2a_5^2a_6a_2  + 59940a_0^2a_5a_6^2a_1 + 330a_0a_5^2a_3^2a_2     + 162a_0a_6a_3^4 - 900a_0^2a_5^2a_4^2 \\
&- 198a_1a_3^3a_2a_6 \\
J_{10} = & \; a_6^{-1}   Res_X \left( f, \frac {\partial f} {\partial X} \right)\\
\end{split}
\end{equation*} 
\begin{equation}\label{eq_J2}
%\begin{multline*}
\begin{split}
& t_2^6 t_3 - 15265260 t_2^5 t_3 - 27949860 t_2^4 t_3^2 - 118098 t_2^3 t_3^3 + 14693280768 t_2^5 + 93437786558880 t_2^4 t_3 \\
& - 878290475269680 t_2^3 t_3^2 + 85811055510240 t_2^2 t_3^3 - 1139016237660 t_2 t_3^4 + 3486784401 t_3^5 \\
& - 223154201664000000 t_2^4  - 287728673929542000000 t_2^3 t_3  - 2469658010168691000000 t_2^2 t_3^2 \\
& - 109818018101695500000 t_2 t_3^3 - 70607384120250000 t_3^4  + 1355661775108800000000000 t_2^3  \\
& + 433843541357670112500000000 t_2^2 t_3 - 662569101476807962500000000 t_2 t_3^2  \\
& + 571919811374025000000000 t_3^3  - 4117822641892980000000000000000 t_2^2  \\
& - 327077365625983809843750000000000 t_2 t_3 - 2316275236064801250000000000000 t_3^2  \\
& + 6253943137374963375000000000000000000 t_2  + 4690457353031222531250000000000000000 t_3   \\
&- 3799270455955290250312500000000000000000000 =  0 \\
\end{split}
%\end{multline*}
\end{equation}
\begin{align*}
  a_{111}  & =  \frac 2 9  \;   (A^2 C - 6 B C + 9 D),\\
  a_{112}  & =  \frac 1 9 \; (2B^3+4ABC+12C^2+3AD),\\
 a_{113}  & =   \,  a_{122} \, = \frac 1 9 \, \left(AB^3 + \frac 4 3 \; A^2 BC + 4 B^2 C + 6 AC^2 + 3 BD \right),\\
   a_{123}  & =  \frac 1 {18}  \, \left(2 B^4 + 4 AB^2C + \frac 4 3 \, A^2 C^2 + 4 BC^2+ 3 ABD+ 12 CD \right),\\
   a_{133}  & =  \frac 1 {18}        \left(AB^4 + \frac 4 3 \, A^2B^2C+ \frac {16} 3 \, B^3 C  + \frac {26} 3 \, ABC^2+ 8 C^3+ 3 B^2 D + 2 ACD\right),\\
   a_{222} & = \frac 1 9    \left( 3B^4 + 6 AB^2 C + \frac 8 3 \, A^2 C^2 + 2 BC^2- 3CD \right),\\
   a_{223} & = \frac 1 {18}    \left(- \frac 2 3 \, B^3 C- \frac 4 3\, ABC^2- 4C^3 + 9 B^2D + 8 ACD \right) ,\\
   a_{233} & = \frac 1 {18}    \left( B^5 + 2AB^3C + \frac 8 9 \, A^2 B C^2 + \frac 2 3 \, B^2 C^2 - BCD + 9 D^2 \right),\\
   a_{333} & = \frac 1 {36} \left( -2 B^4C - 4 A B^2 C^2- \frac {16} 9 \, A^2 C^3 - \frac 4 3 \, BC^3  + 9 B^3 D + 12 ABCD + 20 C^2 D \right) 
\end{align*}
\endgroup

%**********************************************
%\nocite{*}
\bibliographystyle{amsplain}

\bibliography{ref}{}

\end{document}